\documentclass[10pt]{amsart}
\usepackage[psamsfonts]{amssymb}
\usepackage{amsthm,amscd}
\usepackage{type1cm}
\usepackage{color}

\newtheorem{thm}{Theorem}[section]

\newtheorem{lem}{Lemma}[section]
\newtheorem{cor}{Corollary}[section]
\newtheorem{pro}{Proposition}[section]
\theoremstyle{remark}

\theoremstyle{definition}

\allowdisplaybreaks[4]

\title[A classification of endo-commutative straight algebras]{A classification of 2-dimensional endo-commutative straight  algebras of type ${\rm{II_1}}$}

\author[S.-E. Takahasi]{Sin-Ei Takahasi}
\author[K. Shirayanagi]{Kiyoshi Shirayanagi}
\author[M. Tsukada]{Makoto Tsukada}
\address[S.-E. Takahasi]{Laboratory of Mathematics and Games\\ Katsushika 2-371\\ Funabashi\\ Chiba 273-0032\\ Japan}
\address[K. Shirayanagi, M. Tsukada]{Department of Information Science\\ Toho University\\ Miyama 2-2-1\\ Funabashi\\ Chiba 274-8510\\ Japan}
\email{sin\_ei1@yahoo.co.jp}
\email[K. Shirayanagi (Corresponding author)]{kiyoshi.shirayanagi@is.sci.toho-u.ac.jp}
\email[M. Tsukada]{tsukada@is.sci.toho-u.ac.jp}

\makeatletter
\@namedef{subjclassname@2020}{\textup{2020} Mathematics Subject Classification}
\makeatother


\date{\today}

\begin{document}

\begin{abstract}
  In this paper, we present a complete classification of 2-dimensional endo-commutative straight algebras of type II$_1$ over any field.  An endo-commutative algebra is a non-associative algebra in which the square mapping preserves multiplication. A 2-dimensional straight algebra satisfies the condition that there exists an element $x$ such that $x$ and $x^2$ are linearly independent. The term type II$_1$ denotes a distinguishing characteristic of its structure matrix, which has rank 2.
We provide multiplication tables for these algebras, listing them up to isomorphism.
\end{abstract}

\maketitle
%\vspace{5mm}

\section{Introduction}\label{sec:intro}%section 1

Let $A$ be a non-associative algebra.  We define $A$ as {\it{endo-commutative}}, if the square mapping from $A$ to itself preserves multiplication, that is, $x^2y^2=(xy)^2$ holds for all $x, y\in A$. 
In \cite{BCFG}, the class of endo-commutative algebras over a field of characteristic different from 2 and 3 was denoted by $\Omega$, with a primary focus on studying nilalgebras in $\Omega$. However, our main concern lies in the classification problem of endo-commutative algebras over fields of any characteristic. As an initial step, we provided a complete classification of endo-commutative algebras of dimension 2 over the trivial field $\mathbb{F}_2$ with two elements in \cite{TST1}.

We classify two-dimensional algebras into two categories: {\it{curled}} and {\it{straight}}.  A 2-dimensional algebra is curled if the square of any element $x$ is a scalar multiple of $x$; otherwise, it is straight.  In our previous work \cite{TST2}, we presented a complete classification of 2-dimensional endo-commutative curled algebras  over any non-trivial field.  Related research on curled algebras can be found in \cite{level2,length1}.  In \cite{TST3}, we extended our classification to include 2-dimensional endo-commutative straight algebras of rank 1 over any non-trivial field. Here, the rank of an algebra denotes the rank of its structure matrix, which is defined as the matrix of structure constants with respect to a linear basis.  However, one challenging problem remains; the complete classification of 2-dimensional endo-commutative straight algebras of rank 2 over any field.
2-dimensional endo-commutative straight algebras of rank 2 are categorized into three types: types I, II, and III, based on the number of non-zero components in a column of
the structure matrix. Furthermore, type II is further subdivided into three subtypes: II$_1$, II$_2$, and II$_3$, based on which component within the column is the sole zero. For detailed information, refer to Section~\ref{sec:ec-straight-rank2}.  In our subsequent work \cite{TST4}, we presented a complete classification of 2-dimensional endo-commutative straight algebras of type I over any field.  

In this paper, our primary focus is on type II$_1$. We present a complete classification of 2-dimensional endo-commutative straight algebras of type II$_1$ over any field.
For the classifications of associative algebras of dimension 2 over the real and complex number fields, refer to \cite{onkochishin}.
For other studies related to 2-dimensional algebras, see \cite{moduli,2dim,variety,classification}.
In the process of classifying the following, 
we did not find it necessary to introduce any equivalence relation on $K^*:=K\backslash\{0\}$:
\begin{enumerate}
\item 2-dimensional endo-commutative algebras over $\mathbb F_2$ (\cite{TST1}),
\item 2-dimensional endo-commutative curled algebras over any non-trivial field (\cite{TST2}), and
\item 2-dimensional endo-commutative straight algebras of rank 1 over any non-trivial field (\cite{TST3}).
\end{enumerate}
However, when classifying 2-dimensional endo-commutative straight algebras of
type I with rank 2, as demonstrated in \cite{TST4}, 
we needed to introduce a specific equivalence relation on $K^*$. Furthermore, the classification of 2-dimensional endo-commutative straight algebras of type II$_1$
is closely connected to as many as five distinct equivalence relations on $K^*$. Further details can be found in Section~\ref{sec:equiv}.

Our main theorems can be summarized as follows:  If the characteristic ${\rm{char}}\, K$ of the base field $K$ is different from 2, then 2-dimensional endo-commutative straight algebras of type II$_1$ over $K$ are classified into four families:
\[
S_1, S_2, \left\{S_3(t, \varepsilon,\delta)\right\}_{t\in\mathcal H_1, \varepsilon, \delta=\pm1}\, \, {\rm{and}}\, \, \{S_4(t))\}_{t\in K^* \backslash \{\pm1\}}
\]
up to isomorphism.  If ${\rm{char}}\, K=2$, then 2-dimensional endo-commutative straight algebras of type II$_1$ over $K$ are classified into four families:
\[
\{S'_1(t)\}_{t\in\mathcal H_2}, \{S'_2(t)\}_{t\in\mathcal H_3}, \{S'_3(t)\}_{ t\in \mathcal H_4},\, \, {\rm{and}}\, \, \{S'_4(t)\}_{t\in K^*\backslash \{1\}}
\]
up to isomorphism.  Here, $S_i$ and $S'_j$ represent the types of structure matrices. The sets $\mathcal H_1, \mathcal H_2, \mathcal H_3$ and $\mathcal H_4$ denote complete representative systems of the equivalence classes in $K^*$, determined by the four equivalence relations on $K^*$.  Further details are described in
Theorems~\ref{thm:char-not=2} and \ref{thm:char=2}.  As an application of the main theorems, we
also provide complete classifications of 2-dimensional endo-commutative straight algebras of type II$_1$ over specific fields.  Additional details can be found in Corollaries~\ref{cor:sq-rootable}, \ref{cor:K=R}, and \ref{cor:quad-non-extendable}.

The remaining sections are devoted to proving the main theorems. We further partition the type II$_1$ class into three subclasses as defined in Section~\ref{sec:II_1}.
In Sections~\ref{sec:ECS1}, \ref{sec:ECS3}, and
\ref{sec:ECS4}, we classify the algebras within each subclass up to isomorphism, respectively. Each section presents the results for both cases: when ${\rm{char}}\, K$
differs from 2 and when it is equal to 2. To complete the proofs of the main theorems, in Sections~\ref{sec:proof1} and \ref{sec:proof2},
we explore isomorphism relationships between the algebras within any pair of the three subclasses. This exploration is conducted separately
for the cases where ${\rm{char}}\, K$ is not 2 and where ${\rm{char}}\, K$ is 2, respectively.

Throughout this paper, we consider $K$ as an arbitrary field. 

\section{Isomorphism criterion for 2-dimensional algebras}\label{sec:iso-criterion}  %section 2

For any $X=\small{\begin{pmatrix}a&b\\c&d\end{pmatrix}}\in GL_2(K)$, define
\[
\widetilde{X}=\begin{pmatrix}a^2&b^2&ab&ab\\c^2&d^2&cd&cd\\ac&bd&ad&bc\\ac&bd&bc&ad\end{pmatrix},
\]
where $GL_n(K)$ is the general linear group consisting of nonsingular $n\times n$ matrices over $K$.  Then we see that the mapping $X\mapsto\widetilde{X}$ is a group homomorphism from $GL_2(K)$ to $GL_4(K)$ (see \cite{TST2}).

Let $A$ be a 2-dimensional algebra over $K$ with a linear basis $\{e, f\}$.  Then $A$ is determined by the multiplication table $\begin{pmatrix}e^2&ef\\fe&f^2\end{pmatrix}$.  We write
\[
\left\{
    \begin{array}{@{\,}lll}
     e^2=a_1e+b_1f, \\
     f^2=a_2e+b_2f, \\
     ef=a_3e+b_3f, \\
     fe=a_4e+b_4f\\
   \end{array}
  \right. 
\]
with $a_i, b_i\in K\, \, (1\le i\le4)$ and the matrix
$\begin{pmatrix}a_1&b_1\\a_2&b_2\\a_3&b_3\\a_4&b_4\end{pmatrix}$ is called the structure matrix of $A$ with respect to the basis $\{e, f\}$. 

We hereafter will freely use the same symbol $A$ for the matrix and for the algebra because the algebra $A$ is determined by its structure matrix.  For 2-dimensional algebras $A$ and $A'$ over $K$, we see that $A$ and $A'$ are isomorphic iff there is $X\in GL_2(K)$ such that 
\begin{equation}\label{eq:iso}% (1)
A'=\widetilde{X^{-1}}AX.
\end{equation}
This result can be found in \cite{TST2}, and shows the isomorphism criterion for 2-dimensional algebras.

When (\ref{eq:iso}) holds, we say that the matrices $A$ and $A'$ are equivalent and refer to $X$ as a {\it transformation matrix} for the equivalence $A\cong A'$.  Also, we call this $X$ a transformation matrix for the isomorphism $A\cong A'$ as well.
\vspace{3mm}

\section{Characterization of 2-dimensional endo-commutative straight algebras}\label{sec:ec-straight} %section 3

Recall that a 2-dimensional algebra is curled if the square of any element $x$ is a scalar multiplication of $x$, otherwise it is straight.  For a point $(p, q, a, b, c, d)\in K^6$, let $S(p, q, a, b, c, d)$ be an algebra over $K$ with multiplication table
\[
\begin{pmatrix}f&ae+bf\\ce+df&pe+qf\end{pmatrix}
\]
with respect to the linear basis $\{e, f\}$.  The structure matrix of this algebra is given by 
\[
\begin{pmatrix}0&1\\p&q\\a&b\\c&d\end{pmatrix}.
\]
Of course, this algebra is straight.  Conversely, we see easily  that any 2-dimensional straight algebra over $K$ is isomorphic to some $S(p, q, a, b, c, d)$ by replacing the bases.  The following results can be found in \cite{TST3}.
\vspace{2mm}

\begin{lem}\label{lem:ec}% Lemma 3.1
The straight algebra $S(p, q, a, b, c, d)$ is endo-commutative iff the point $(p, q, a, b, c, d)\in K^6$ satisfies 
\begin{equation}\label{eq:ec}% (2)
\left\{\begin{array}{@{\,}lll}   
pq+pc=pb^2+a^2b+abc, \\
p(c-a)=(b-d)\{p(b+d)-q(a+c)\}, \\
p(d-b)=a^2-c^2, \\
q^2+pd=a^2+qb^2+ab^2+abd, \\
q(d-b)=ab-cd.\\  
\end{array} \right. 
\end{equation}
\end{lem}

\begin{lem}\label{lem:iso}% Lemma 3. 2
The straight algebras $S(p, q, a, b, c, d)$ and $S(p', q', a', b', c', d')$ are isomorphic iff there are $x, y, z, w\in K$ with $\left|\begin{matrix}x&y\\z&w\end{matrix}\right|\ne0$ such that 
\begin{equation}\label{eq:st-iso}% (3)
\left\{\begin{array}{@{\,}lll}
p'y^2+(a'+c')xy=z,\\
x^2+q'y^2+(b'+d')xy=w,\\
p'w^2+(a'+c')zw=px+qz, \\
z^2+q'w^2+(b'+d')zw=py+qw,\\
p'yw+a'xw+c'yz=ax+bz, \\
xz+q'yw+b'xw+d'yz=ay+bw,\\
p'yw+a'yz+c'xw=cx+dz,\\
xz+q'yw+b'yz+d'xw=cy+dw
\end{array} \right. 
\end{equation}
holds.
\end{lem}

Note that in Lemma~\ref{lem:iso}, $X:=\begin{pmatrix}x&y\\z&w\end{pmatrix}$ is a transformation matrix  for the isomorphism $S(p, q, a, b, c, d)\cong S(p', q', a', b', c', d')$.

\section{2-dimensional endo-commutative straight algebras of rank 2}\label{sec:ec-straight-rank2} %section 4

Let $\mathcal{EC}$ be the family of all 2-dimensional endo-commutative algebras over $K$ with a linear basis $\{e, f\}$ and put
\[
\mathcal{ECS}_{{\rm{(rank \, 2)}}}=\{S(p, q, a, b, c, d)\in\mathcal{EC} : {\rm{rank}}\, S(p, q, a, b, c, d)=2\}.
\]
Define
\[
\mathcal{ECS}_{{\rm{I}}}=\mathcal{ECS}_{001}\sqcup\mathcal{ECS}_{010}\sqcup\mathcal{ECS}_{100},
\]
where 
\[
\left\{\begin{array}{@{\,}lll}  
\mathcal{ECS}_{001}=\{S(p, q, a, b, c, d)\in\mathcal{EC} : p=a=0, c\ne0\}, \\
\mathcal{ECS}_{010}=\{S(p, q, a, b, c, d)\in\mathcal{EC} : p=c=0, a\ne0\}, \\
\mathcal{ECS}_{100}=\{S(p, q, a, b, c, d)\in\mathcal{EC} : a=c=0, p\ne0\}.
\end{array} \right. 
\]
Then $\mathcal{ECS}_{{\rm{I}}}$ is a subfamily of $\mathcal{ECS}_{{\rm{(rank \, 2)}}}$.  We say that each algebra in $\mathcal{ECS}_{{\rm{I}}}$ is of type I.   Define
\[
\mathcal{ECS}_{{\rm{II}}}=\mathcal{ECS}_{011}\sqcup\mathcal{ECS}_{101}\sqcup\mathcal{ECS}_{110},
\]
where 
\[
\left\{\begin{array}{@{\,}lll}  
\mathcal{ECS}_{011}=\{S(p, q, a, b, c, d)\in\mathcal{EC} : p=0, a, c\ne0\}, \\
\mathcal{ECS}_{101}=\{S(p, q, a, b, c, d)\in\mathcal{EC} : a=0, p, c\ne0\}, \\
\mathcal{ECS}_{110}=\{S(p, q, a, b, c, d)\in\mathcal{EC} : c=0, p, a\ne0\}.
\end{array} \right. 
\]
Then $\mathcal{ECS}_{{\rm{II}}}$ is a subfamily of $\mathcal{ECS}_{{\rm{(rank \, 2)}}}$.  We say that each algebra in $\mathcal{ECS}_{{\rm{II}}}$ is of type II.
Furthermore, we say that each algebra in $\mathcal{ECS}_{011}, \mathcal{ECS}_{101}$ and $\mathcal{ECS}_{110}$ is of type ${\rm{II_1}}, {\rm{II_2}}$ and ${\rm{II_3}}$, respectively.  Define
\[
\mathcal{ECS}_{{\rm{III}}}=\{S(p, q, a, b, c, d)\in\mathcal{EC} : p, a, c\ne0\}.
\]
Then $\mathcal{ECS}_{{\rm{III}}}$ is a subfamily of $\mathcal{ECS}_{{\rm{(rank \, 2)}}}$.  We say that each algebra in $\mathcal{ECS}_{{\rm{III}}}$ is of type III.  Since the algebra $S(p, q, a, b, c, d)$ has rank 2 iff at least one of  $p, a$ and $c$ is nonzero, we see easily that
\[
\mathcal{ECS}_{{\rm{(rank \, 2)}}}=\mathcal{ECS}_{{\rm{I}}}\sqcup\mathcal{ECS}_{{\rm{II}}}\sqcup\mathcal{ECS}_{{\rm{III}}}.
\]
In this paper, we focus on the classification problem of the type II$_1$ algebras, that is, $\mathcal{ECS}_{011}$.
\vspace{1mm}

\section{Characterization of type II$_1$ algebras}\label{sec:II_1} %section 5

Recall that
\[
\mathcal{ECS}_{011}=\{S(p, q, a, b, c, d)\in\mathcal{EC} : p=0, a, c\ne0\}.
\]

Let $S(p, q, a, b, c, d)\in\mathcal {ECS}_{011}$.  Then $p=0, a, c\ne0$ and the point $(p, q, a, b, c, d)\in K^6$ satisfies (\ref{eq:ec}) by Lemma~\ref{lem:ec}.  Since $p=0$ and $a, c\ne0$, (\ref{eq:ec}) is rewritten as
\begin{equation}\label{eq:II_1}% (4)
\left\{\begin{array}{@{\,}lll}   
b(a+c)=0, \\
qd(a+c)=0, \\
a^2=c^2, \\
q^2=a^2+qb^2+ab^2+abd, \\
q(d-b)=ab-cd.\\  
\end{array} \right. 
\end{equation}
Define
\[
\left\{\begin{array}{@{\,}lll}  
(\mathcal{ECS}_{011})_1=\{S(p, q, a, b, c, d)\in\mathcal{ECS}_{011} : b=0\}\\
(\mathcal{ECS}_{011})_{2}=\{S(p, q, a, b, c, d)\in\mathcal{ECS}_{011} : b\ne0, q=0\}\\
(\mathcal{ECS}_{011})_{3}=\{S(p, q, a, b, c, d)\in\mathcal{ECS}_{011} : b\ne0, q\ne 0, d=0\}\\
(\mathcal{ECS}_{011})_{4}=\{S(p, q, a, b, c, d)\in\mathcal{ECS}_{011} : b\ne0, q\ne0, d\ne0\}.
\end{array} \right. 
\]
Then we have
\[
\mathcal{ECS}_{011}=(\mathcal{ECS}_{011})_{1}\sqcup(\mathcal{ECS}_{011})_{2}\sqcup(\mathcal{ECS}_{011})_{3}\sqcup(\mathcal{ECS}_{011})_{4}.
\]

\begin{lem}\label{lem:II_1_1} % Lemma 5. 1
\[
(\mathcal{ECS}_{011})_{1}=\{S(0, a, a, 0, -a, d) :  a, d\ne0\}\sqcup\{S(0, \varepsilon a, a, 0, \delta a, 0) : a\ne0, \varepsilon, \delta=\pm1\}.
\]
\end{lem}

\begin{proof}
Let $(p, q, a, b, c, d)\in K^6$.  Suppose (\ref{eq:II_1}) and $p=b=0, a, c\ne0$.  Then (\ref{eq:II_1}) is rewritten as
\begin{equation}\label{eq:II_1_1}% (5)
\left\{\begin{array}{@{\,}lll}   
qd(a+c)=0, \\
a^2=c^2, \\
q^2=a^2, \\
qd=-cd. 
\end{array} \right. 
\end{equation}
Since $a\ne0$, the third equation of (\ref{eq:II_1_1}) implies $q\ne0$, and hence (\ref{eq:II_1_1}) is rewritten as
\begin{equation}\label{eq:II_1_1a} % (6)
\left\{\begin{array}{@{\,}lll}   
d(a+c)=0, \\
a^2=c^2, \\
q^2=a^2, \\
d(q+c)=0. 
\end{array} \right. 
\end{equation}
If $d\ne0$, then (\ref{eq:II_1_1a}) is rewritten as $\left\{\begin{array}{@{\,}lll}c=-a, \\q=a.\end{array} \right.$  If $d=0$, then (\ref{eq:II_1_1a}) is rewritten as $\left\{\begin{array}{@{\,}lll}c=\delta a, \\q=\varepsilon a,\end{array} \right.$ where $\varepsilon, \delta=\pm1$.  Then we obtain the desired result.
\end{proof}
 
\begin{lem}\label{lem:II_1_2} % Lemma 5. 2
$(\mathcal{ECS}_{011})_{2}=\emptyset$.
\end{lem}

\begin{proof}
Suppose $(\mathcal{ECS}_{011})_{2}\ne\emptyset$ and take $S(p, q, a, b, c, d)\in(\mathcal{ECS}_{011})_{2}$.  Then the point $(p, q, a, b, c, d)\in K^6$ satisfies (\ref{eq:II_1}).  Since $b\ne0$ and $q=0$, (\ref{eq:II_1}) is rewritten as
\begin{equation}\label{eq:II_1_2} % (7)
\left\{\begin{array}{@{\,}lll}   
a+c=0, \cdots(7-1)\\
0=a^2+ab^2+abd, \cdots(7-2)\\
0=ab-cd.\cdots(7-3)\\  
\end{array} \right. 
\end{equation}
By (7-1) and (7-3), we have $a(b+d)=0$, hence $d=-b$ because $a\ne0$. Substituting this into (7-2), we have $a^2=0$, a contradiction.
\end{proof}

\begin{lem}\label{lem:II_1_3} % Lemma 5. 3
$(\mathcal{ECS}_{011})_{3}=\{S(0, -a, a, b, -a, 0) : a, b\ne0\}$.
\end{lem}

\begin{proof}
Let $(p, q, a, b, c, d)\in K^6$.  Suppose (\ref{eq:II_1}) and $p=d=0, a, c, b, q\ne0$.  Then  (\ref{eq:II_1}) is rewritten as
\[
\left\{\begin{array}{@{\,}lll}   
a+c=0, \\
q^2=a^2+qb^2+ab^2, \\
-q=a, 
\end{array} \right. 
\]
which is equivalent to 
$\left\{\begin{array}{@{\,}lll}   
c=-a, \\
q=-a.
\end{array} \right.$  Then we obtain the desired result.
\end{proof}

\begin{lem}\label{lem:II_1_4} % \Lemma 5. 4
The following two statements hold.
\vspace{2mm}

{\rm{(i)}} Suppose ${\rm{char}}K\ne2$.  Then 
\begin{align*}
(\mathcal{ECS}_{011})_{4}&=\{S(0, (b+d)^2/4, (d^2-b^2)/4, b, (b^2-d^2)/4, d) : b, d\ne0, b\ne\pm d\}.
\end{align*}

{\rm{(ii)}} Suppose that ${\rm{char}}K=2$.  Then 
\[
(\mathcal{ECS}_{011})_{4}=\{S(0, q, a, b, a, b) : q, a, b\ne0, q^2+a^2+qb^2=0\}.
\]
\end{lem}

\begin{proof}
Let $S(p, q, a, b, c, d)\in(\mathcal{ECS}_{011})_{4}$.  Then the point $(p, q, a, b, c, d)\in K^6$ satisfies (\ref{eq:II_1}).  Since $b, q, d\ne0$, (\ref{eq:II_1}) is rewritten as
\begin{equation}\label{eq:II_1_4} %(8)
\left\{\begin{array}{@{\,}lll}   
a+c=0, \cdots(8-1)\\
q^2=a^2+qb^2+ab^2+abd, \cdots(8-2)\\
q(d-b)=ab-cd.\cdots (8-3)\\  
\end{array} \right. 
\end{equation}

(A) The case where $b=d$.  By (8-1), we see that (8-3) becomes $a(b+d)=0$, hence $b+d=0$ (because $a\ne0$), so $2b=0$.  If ${\rm{char}}K\ne2$, then $b=0$, a contradiction. If ${\rm{char}}K=2$, then (\ref{eq:II_1_4}) becomes
$\left\{\begin{array}{@{\,}lll}   
a=c, \\
q^2=a^2+qb^2.
\end{array} \right. $

(B) The case where $b\ne d$.  By (8-1) and (8-3), we get $q=a(b+d)/(d-b)$.  Substituting this into (8-2), we get
\[
\frac{a^2(b+d)^2}{(d-b)^2}=a^2+\frac{a(b+d)b^2}{d-b}+ab^2+abd.
\]
Reducing $a$ (since $a\ne0$) and clearing the denominators, the above equation becomes $4abd=bd^3-b^3d$, hence $4a=d^2-b^2$ because $bd\ne0$.  Consequently, (\ref{eq:II_1_4}) is rewritten as 
\[
\left\{\begin{array}{@{\,}lll}   
c=-a, \\
4a=d^2-b^2, \\
q=a(b+d)/(d-b).
\end{array} \right. 
\]
If ${\rm{char}}K\ne2$, then (\ref{eq:II_1_4}) is rewritten as 
\[
\left\{\begin{array}{@{\,}lll}   
c=-a, \\
a=(d^2-b^2)/4, \\
q=(b+d)^2/4.
\end{array} \right. 
\]
Also if ${\rm{char}}K=2$, then $4a=d^2-b^2$ implies $(d-b)^2=0$, hence $b=d$, a contradiction. 

Summarizing (A) and (B) as presented above, we achieve the desired result.  
\end{proof}

From now on let $K^*=K\, \backslash\, \{0\}$, which is a multiplicative group.
\vspace{2mm}

\section{Five equivalence relations on $K^*$}\label{sec:equiv} %section 6
\vspace{2mm}

As seen in \cite{TST4}, the classification problem of the type I algebras is closely tied to a specific equivalence relation on $K^*$. This holds true for the classification problem of the type II$_1$ algebras as well.  In this section, we will enumerate the necessary equivalence relations and provide thorough analysis for each of them.

Define 
\[
(K^*)^2:=\{k^2 : k\in K^*\}\, \, {\rm{and}}\, \, K_{\iota^2+\iota}:=\{x^2+x : x\in K\}.
\]

(I) If $t, t'\in K^*$ satisfy $t/t'\in(K^*)^2$, we write $t\sim_1t'$.  Since $(K^*)^2$ is a group, it is obvious that $\sim_1$ is an equivalence relation on $K^*$.  Let $\mathcal H_1$ denote a complete representative system of the quotient set $K^*/\sim_1$ of $K^*$ by $\sim_1$.  

We refer to $K$ as square-rootable if $K=\{k^2 : k\in K\}$.  For example, the complex number field $\mathbb C$ is square-rootable.  However, the rational number field $\mathbb Q$ and the real number field $\mathbb R$ are not square-rootable.  
\vspace{2mm}

\begin{lem}\label{lem:equiv1}% Lemma 6. 1
{\rm{(i)}} If $K$ is square-rootable, we can let $\mathcal H_1=\{1\}$.
\vspace{2mm}

{\rm{(ii)}} If $K=\mathbb R$, we can let $\mathcal H_1=\{-1, 1\}$.
\vspace{2mm}

{\rm{(iii)}} If $K=\mathbb Q$, then $\, \sharp(K^*/\sim_{1})=\aleph_0$, hence any complete representative system of $K^*/\sim_{1}$ is countably infinite.
\end{lem}

\begin{proof}
(i)  Suppose $K$ is square-rootable.  Then, for any $t'\in K^*$, we have $1/t'\in K^*=(K^*)^2$, hence $1\sim_1t'$.
\vspace{2mm}

(ii) Suppose $K=\mathbb R$ and let $t, t'\in K^*$.  Then $t\sim_1t'$ iff $t/t'>0$.  Therefore we see that $1\sim_1t\,  (\forall t>0)$, $-1\sim_1t\,  (\forall t<0)$ and $1\nsim_1-1$.  This observation implies that $\{-1, 1\}$ is a complete representative system of $K^*/\sim_1$.  
\vspace{2mm}

(iii)  Suppose $K=\mathbb Q$ and let $\mathbb P$ be the set of all prime numbers, hence $\mathbb P\subset \mathbb Q^*$.  To show $\, \sharp(K^*/\sim_{1})=\aleph_0$, it suffices to show that $p\nsim_1q$ holds for all $p, q\in\mathbb P$ with $p\ne q$.  We assume that $p\sim_1q$ for some $p, q\in\mathbb P$ with $p\ne q$, and derive the contradiction as follows.  By the definition of $\sim_1$, there exist $m, n\in\mathbb N$ such that $pm^2=qn^2$.
However, when performing prime factorization on $pm^2$, we find that the prime factor $p$ appears an odd number of times in total. On the other hand, when factorizing
$qn^2$, if the prime factor $p$ is present, it appears an even number of times in total. This is a contradiction.
\end{proof}

(II) Suppose ${\rm{char}}\, K=2$.  If $t, t'\in K^*$ satisfy $t+t'\in K_{\iota^2+\iota}$, we write $t\sim_{2}t'$.  Then $\sim_{2}$ is an equivalence relation on $K^*$.  Indeed, both reflexivity and symmetry trivially hold.  To show the transitivity, suppose that $t\sim_{2}t'$ and $t'\sim_{2}t''$.  Then there are $x, y\in K$ such that $ t+t'=x^2+x$ and $t'+t''=y^2+y$. Therefore
\[
t+t''=t+t'+t'+t''=x^2+x+y^2+y=(x+y)^2+x+y\in K_{\iota^2+\iota},
\]
hence $t\sim_{2}t''$, as required.   Let $\mathcal H_2$ denote a complete representative system of the  quotient set $K^*/\sim_{2}$.

If any quadratic polynomial over $K$ has a root in $K$, then $K$ is said to be quadratically non-extendable.  Otherwise, $K$ is said to be quadratically extendable. 

If $K$ is quadratically non-extendable, then we can obviously let $\mathcal H_2=\{1\}$.

Let $\mathbb F_2[X]$ be the ring of all polynomials over the trivial field $\mathbb F_2$ of two elements  and $\mathbb F_2(X)$ the filed of all rational functions over $\mathbb F_2$.  Then $\mathbb F_2(X)$ is obviously an infinite  field of characteristic 2.  Also, $\mathbb F_2(X)$ is quadratically extendable.  In fact, we have only to show that the quadratic polynomial $x^2+x+1$ does not have any roots in $\mathbb F_2(X)$.  Otherwise, there exist $p(X)$ and $q(X)\in\mathbb F_2[X]\backslash\{0\}$ such that
$1=(p(X)/q(X))^2+p(X)/q(X)$, that is, $q^2=(p+q)p$.  For any $f(X)\in\mathbb F_2[X]$, we denote by $d(f)$ the degree of $f(X)$.  If $d(p)=d(q)$, then $d((p+q)p)<2d(p)$ because the leading terms of $p(X)$ and $q(X)$ are both $X^{d(p)}$ and so their sum is 0.  Then $2d(q)=d(q^2)=d((p+q)p)<2d(p)$, a contradiction.  If $d(p)\ne d(q)$,
then $2d(q)=d(q^2)=d((p+q)p)=\max\{d(p), d(q)\}+d(p)$, hence $d(q)=d(p)$, a contradiction. Moreover, we have the following.

\begin{lem}\label{lem:equiv2}% Lemma 6. 2
If $K=\mathbb F_2(X)$, then $\, \sharp(K^*/\sim_{2})=\aleph_0$, hence any complete representative system of $K^*/\sim_{2}$ is countably infinite.
\end{lem}

\begin{proof}
It suffices to show that for any $i\in\mathbb N$, $X^{2i+1}$ can be taken as a representative of an equivalence class in $\mathbb F_2(X)^*/\sim_{2}$, so that distinct $i$'s give distinct representatives.  For this, we will prove that if $i\ne j$, then $X^{2i+1}\nsim_{2}X^{2j+1}$.  We assume that for some $i$ and $j$ with $i<j$, $X^{2i+1}\sim_{2}X^{2j+1}$ and derive the contradiction as follows.  By the definition of $\sim_2$, there exist $p(X), q(X)\in\mathbb F_2[X]\backslash\{0\}$ such that 
\[
X^{2i+1}+X^{2j+1}=\left(\frac{p(X)}{q(X)}\right)^2+\frac{p(X)}{q(X)}.
\]
By the above equality, we have
\[
X^{2i+1}q(X)^2+X^{2j+1}q(X)^2=p(X)^2+p(X)q(X),
\]
that is, $X^{2i+1}q^2+X^{2j+1}q^2=p^2+pq$.  Put $m=d(p), n=d(q)$, and $D=d(p^2+pq)$.  Since $d(X^{2i+1}q^2+X^{2j+1}q^2)=\max\{2i+1, 2j+1\}+2n=(2j+1)+2n$, it follows that $(2j+1)+2n=D$.
\vspace{2mm} 

$\bullet$  If $m>n$, then $D=2m$, so $(2j+1)+2n=2m$ and hence $1\in2\mathbb Z$, which is a contradiction.
\vspace{2mm} 
 
$\bullet$ If $m=n$, then $D<2n$ since the leading terms of $p(X)^2$ and $p(X)q(X)$ are both $X^{2n}$ and so their sum is 0.  Therefore, $(2j+1)+2n<2n$, which is a contradiction.
\vspace{0mm} 
 
$\bullet$ If $m<n$, then $D=m+n$, so $(2j+1)+2n=m+n$ and hence $m=n+(2j+1)>n$, which is a contradiction.
\end{proof}

(III) Suppose ${\rm{char}}\, K=2$.  Let $t, t'\in K^*$.  If there are $x\in K^*$ and $y\in K$ such that $t'x^2+y^2+t=0$, we write $t\sim_3t'$.  Then $\sim_3$ is an equivalence relation on $K^*$.  Indeed, if $t\in K^*$, then $tx^2+y^2+t=2t=0$ for $x=1$ and $y=0$, hence $t\sim_3t$, that is,  reflexivity follows.  To show symmetry, suppose $t, t'\in K^*$ with $t\sim_3t'$.  Then there are $x\in K^*$ and $y\in K$ such that $t'x^2+y^2+t=0$.  Define $\left\{\begin{array}{@{\,}lll}X=1/x, \\Y=y/x,\end{array} \right. $ hence $X\in K^*$ and $Y\in K$.  Also, we have 
\begin{align*}
tX^2+Y^2+t'&=t/x^2+y^2/x^2+t'=(t'x^2+y^2)/x^2+y^2/x^2+t'\\
&=2t'+2y^2/x^2=0,
\end{align*}
hence $t'\sim_3t$, that is, symmetry follows.  To demonstrate transitivity, suppose that $t\sim_3t'$ and $t'\sim_3t''$.  Then there are $x, x'\in K^*$ and $y, y'\in K$ such that $t'x^2+y^2+t=0$ and $t''x'^2+y'^2+t'=0$.   Define $\left\{\begin{array}{@{\,}lll}X=x'x, \\Y=y'x+y,\end{array} \right. $ hence $X\in K^*$ and $Y\in K$.  Also, we have  
\begin{align*}
t''X^2+Y^2+t&=t''x'^2x^2+(y'x+y)^2+t=t''x'^2x^2+y'^2x^2+y^2+t\\
&=(y'^2+t')x^2+y'^2x^2+y^2+t=t'x^2+y^2+t=0,
\end{align*}
hence $t\sim_3t''$, that is, transitivity follows.   Let $\mathcal H_3$ denote a complete representative system of  $K^*/\sim_3$.  

If $K$ is square-rootable, then we can obviously let $\mathcal H_3=\{1\}$.   Moreover, we have the following.
 
 \begin{lem}\label{lem:equiv3}% Lemma 6.3
 If $K=\mathbb F_2(X)$, then $K^*/\sim_3\, =\{\widetilde{1}, \widetilde{X}\}$ and $\widetilde{1}\ne\widetilde{X}$, where $\widetilde{f}$ denotes the coset of $f\in K^*$.
 \end{lem}

\begin{proof}
First of all, we show that the following statements.
\vspace{2mm}

(i) If $p(X)\in\mathbb F_2[X]$ has a term of odd degree, then $p(X)\sim_3X$.
\vspace{2mm}

(ii) If $p(X)\in\mathbb F_2[X]$ has no term of odd degree, then $p(X)\sim_31$.
\vspace{2mm}

In fact, let $p(X)\in\mathbb F_2[X]$.  To show (i), suppose that $p(X)$ has a term of odd degree.  By splitting $p(X)$ into the terms of odd degree and the terms of even degree, we have
\[
p(X)=X^{2i_1+1}+\cdots+X^{2i_k+1}+q(X),
\]
where $q(X)=0$ if $p(X)$ has no term of even degree, and $q(X)=X^{2j_1}+\cdots+X^{2j_l}$ if $p(X)$ has a term of even degree.  Then we see that $q(X)=r(X)^2$ for some $r(X)\in\mathbb F_2[X]$, hence
\[
X\cdot(X^{i_1}+\cdots+X^{i_k})^2+r(X)^2+p(X)=0,
\]
so $p(X)\sim_3X$ .  Thus (i) follows.  Next, to show (ii), suppose that $p(X)\in\mathbb F_2[X]$ has no term of odd degree, and hence we can write 
\[
p(X)=X^{2j_1}+\cdots+X^{2j_l}.
\]
 Then we have
 \[
 1\cdot(X^{j_1}+\cdots+X^{j_l})^2+0^2+p(X)=0,
 \]
hence $p(X)\sim_31$.  Thus (ii) follows.  Next we show the following.
\vspace{2mm}

(iii) $p(X)/q(X)\sim_3p(X)q(X)$ holds for all $p(X), q(X)\in\mathbb F_2[X]$ with $q(X)\ne0$.
\vspace{2mm}

In fact, let  $p(X), q(X)\in\mathbb F_2[X]$ with $q(X)\ne0$.  Since
\[
\left(\frac{p(X)}{q(X)}\right)\cdot q(X)^2+0^2+p(X)q(X)=0, 
\]
it follows that  $p(X)/q(X)\sim_3p(X)q(X)$.  Thus (iii) follows.

Now, by (i), (ii) and (iii), we can easily see that any non-zero rational function $h(X)$ is equivalent to either $1$ or $X$, that is, $h(X)\sim_31$ or $h(X)\sim_3 X$.

Finally, to complete the proof, we will show $X\nsim_31$.  In fact, assume $X\sim_31$.  Then there is $t\in K^*$ and $s\in K$ such that $1\cdot t^2+s^2+X=0$.  Then $t+s=p(X)/q(X)$ for some $p(X), q(X)\in\mathbb F_2[X]$ with $q(X)\ne0$, hence 
\[
Xq(X)^2=(t+s)^2q(X)^2=\left(\frac{p(X)}{q(X)}\right)^2q(X)^2=p(X)^2.
\]
The degree of the polynomial on the leftmost side is odd, whereas the degree of 
the polynomial on the rightmost side is even. This is a contradiction. 
\end{proof}

By Lemma~\ref{lem:equiv3}, we see that any complete representative system of $\mathbb F_2(X)^*/\sim_{2}$ consists of two elements.
The following lemma implies that there exists a field of characteristic 2 in which any complete representative system with respect to $\sim_3$ is infinite.

 \begin{lem}\label{lem:equiv3-infinity}% Lemma 6. 4
Let $K:=\mathbb F_2((X_i)_{i\in\mathbb N})$ be the field of rational functions in infinitely many variables over $\mathbb F_2$.  Then $\, \sharp(K^*/\sim_{3})=\aleph_0$, hence any complete representative system of $K^*/\sim_{3}$ is countably infinite.
 \end{lem}

\begin{proof}
It suffices to prove that for any $i, j\in\mathbb N$, if $i\ne j$, then $X_i\nsim_3X_j$.  Assume $X_i\sim_3X_j$.  Then there are $t\in K^*$ and $s\in K$ such that $X_jt^2+s^2+X_i=0$.  Let $t=h/g$ and $s=q/p$, where $g, h, p\in \mathbb F_2[(X_i)_{i\in\mathbb N}]\backslash\{0\}$ and $q\in \mathbb F_2[(X_i)_{i\in\mathbb N}]$.  Thus $X_i=X_j(h/g)^2+(q/p)^2$.  By clearing the denominators, we have
\[
X_ig^2p^2=X_jh^2p^2+g^2q^2.
\]
Hence, the degree of the polynomial on the left hand side with respect to $X_i$ is odd, whereas the degree of the polynomial on the right hand side with respect to $X_i$ is even. This is a contradiction.
\end{proof}

(IV)  Suppose ${\rm{char}}\, K=2$.  If $t, t'\in K^*$ satisfy $1/t+1/t'\in K_{\iota^2+\iota}$, we write $t\sim_4t'$.  Then we see easily that  $\sim_4$ is an equivalence relation on $K^*$ as observed in (II).  Let $\mathcal H_4$ denote a complete representative system of $K^*/\sim_4$.

If $K$ is quadratically non-extendable, then we can obviously let $\mathcal H_4=\{1\}$ as observed in (II).  Also, we have the following.

\begin{lem}\label{lem:equiv5} % Lemma 6. 5

If $K=\mathbb F_2(X)$, then $K$ is quadratically extendable and $\, \sharp(K^*/\sim_{4})=\aleph_0$, hence any complete representative system of $K^*/\sim_{4}$ is countably infinite.
\end{lem}

\begin{proof}
Suppose $K=\mathbb F_2(X)$.  Then $K$ is quadratically extendable as observed in (II).  Next we have only to show that $\{X^{-(2i+1)} : i\in\mathbb N\}$ is a representative system of $K^*/\sim_4$.  To do this, let $i, j\in\mathbb N$ with $i\ne j$.  Then $X^{2i+1}\nsim_2X^{2j+1}$ as shown in the proof of Lemma~\ref{lem:equiv2}.  Then we see easily that $X^{-(2i+1)}\nsim_5X^{-(2j+1)}$, hence  $\{X^{-(2i+1)} : i\in\mathbb N\}$ is a representative system of $K^*/\sim_4$, as required.  
\end{proof}

(V)  Suppose ${\rm{char}}\, K\ne2$.  If $t, t'\in K^*\backslash\{-4\}$ satisfy $\frac{t'(4+t)}{t(4+t')}\in(K^*)^2$, we write $t\sim_5t'$. Then we see easily that  $\sim_5$ is an equivalence relation on $K^*\backslash\{-4\}$ because $(K^*)^2$ is a group.  Let $\mathcal H'_5$ denote a complete representative system of $(K^*\, \backslash\, \{-4\})/\sim_5$. 

When $K$ is square-rootable, we can obviously take $\mathcal H'_5=\{1\}$ as a complete representative system of  $K^*\, \backslash\, \{-4\}$.

\begin{lem}\label{lem:equiv4} % Lemma 6. 6
{\rm(i)}  If $K=\mathbb R$, we can take $\mathcal H'_5=\{-1, 1\}$ as a complete representative system of  $K^*\, \backslash\, \{-4\}$.  
\vspace{2mm}

{\rm(ii)}  If $K=\mathbb Q$, then $\, \sharp(K^*/\sim_{5})=\aleph_0$, hence any complete representative system of  $K^*\, \backslash\, \{-4\}$ is countably infinite.  

\end{lem}

\begin{proof}
(i)  Suppose $K=\mathbb R$ and let $t, t'\in K^*\backslash\{-4\}$.  Then 
\begin{align*}
t\sim_4t'&\Leftrightarrow \frac{t'(4+t)}{t(4+t')}>0\Leftrightarrow
\left\{\begin{array}{@{\,}lll}   
t(4+t)>0\\
t'(4+t')>0
\end{array} \right. 
{\rm{or}}\, \, 
\left\{\begin{array}{@{\,}lll}   
t(4+t)<0\\
t'(4+t')<0
\end{array} \right.\\
&\Leftrightarrow 
\left\{\begin{array}{@{\,}lll}   
t<-4\, \, {\rm{or}}\,\, 0<t\\
t'<-4\,\,  {\rm{or}}\,\,  0<t'
\end{array} \right. 
{\rm{or}}\, \, 
\left\{\begin{array}{@{\,}lll}   
-4<t<0\\
-4<t'<0
\end{array} \right.\\
&\Leftrightarrow  {\rm{Either}}\, \, C_1(t, t'), C_2(t, t'), C_3(t, t'), C_4(t, t')\,\,  {\rm{or}}\,\, C_5(t, t')\, \,  {\rm{holds}},
\end{align*}
where $C_1(t, t')\equiv\left\{\begin{array}{@{\,}lll}0<t\\0<t',\end{array} \right. C_2(t, t')\equiv\left\{\begin{array}{@{\,}lll}t<-4\\t'<-4,\end{array} \right. C_3(t, t')\equiv\left\{\begin{array}{@{\,}lll}0<t\\t'<-4,\end{array} \right. C_4(t, t')\equiv\left\{\begin{array}{@{\,}lll}-4<t<0\\-4<t'<0\end{array} \right.$ and  $C_5(t, t')\equiv\left\{\begin{array}{@{\,}lll}t<-4\\0<t'.\end{array} \right.$
Hence we can take $\mathcal H'_4=\{-1, 1\}$ as a complete representative system of  $K^*\, \backslash\, \{-4\}$.  \vspace{2mm}

(ii)  Suppose $K=\mathbb Q$ and take a sequence $\{p_1, p_2, \cdots\}$ in $\mathbb P$ such that $p_i+4<p_{i+1}$ for all $i\in\mathbb N$.  Then it suffices to show that $p_i\nsim_5p_j$ for all $i, j\in\mathbb N$ with $i<j$.  We assume that $p_i\sim_5p_j$ for some $i, j\in\mathbb N$ with $i<j$, and derive the contradiction as follows.  By hypothesis, there are $m, n\in\mathbb N$ with $p_j(4+p_i)m^2=p_i(4+p_j)n^2$.
However, when performing prime factorization on $p_j(4+p_i)m^2$, the prime factor $p_j$ appears an odd number of times in total. On the other hand,
when factorizing $p_i(4+p_j)n^2$, if the prime factor $p_j$ is present, it appears an even number of times in total. This is a contradiction.
\end{proof}

\section{Main results}\label{sec:main-results} %section 7

Let us retain the symbols $\mathcal H_i\, \, (1\le i\le4) $ as they are in the previous section. Then the main theorems are stated as follows.
\vspace{2mm}

\begin{thm}\label{thm:char-not=2} % Theorem 7. 1
Suppose ${\rm{char}}\, K\ne2$.  Up to isomorphism, 2-dimensional endo-commutative straight algebras of type II$_1$ over $K$ are classified into 4 families
\[
S(0, 1, 1, 0, -1, 2), S(0, 4, -4, -4, 4, 0), \left\{S(0, \varepsilon t, t, 0, \delta t, 0)\right\}_{t\in\mathcal H_1, \varepsilon, \delta=\pm1}
\]
and
\[
 \{S(0, (1+t)^2/4, (t^2-1)/4, 1, (1-t^2)/4, t)\}_{t\in K^* \backslash \{\pm1\}}
\]
with multiplication tables on a linear basis $\{e, f\}$ defined by
\[
\begin{pmatrix}f&e\\-e+2f&f\end{pmatrix}, \begin{pmatrix}f&-4(e+f)\\4e&4f\end{pmatrix}, \begin{pmatrix}f&te\\\delta te&\varepsilon tf\end{pmatrix}\, \,{\rm{and}}\, \,  \begin{pmatrix}f&\frac{t^2-1}{4}e+f\\\frac{1-t^2}{4}e+tf&\frac{(1+t)^2}{4}f\end{pmatrix},
\]
respectively.
\end{thm}
\vspace{1mm}

If $K$ is square-rootable and ${\rm{char}}\, K\ne2$, then $\mathcal H_1=\{1\}$ from Lemma~\ref{lem:equiv1} (i). Then
Theorem~\ref{thm:char-not=2} implies the following.

\begin{cor}\label{cor:sq-rootable} % Corollary 7. 1
Suppose that $K$ is square-rootable and ${\rm{char}}\, K\ne2$.  Up to isomorphism, 2-dimensional endo-commutative straight algebras of type II$_1$ over $K$ are classified into 4 families $S(0, 1, 1, 0, -1, 2)$, $S(0, 4, -4, -4, 4, 0)$, $\left\{S(0, \varepsilon, 1, 0, \delta, 0)\right\}_{\varepsilon, \delta=\pm1}$ and $ \{S(0, (1+t)^2/4, (t^2-1)/4, 1, (1-t^2)/4, t)\}_{t\in K^* \backslash \{\pm1\}}$.
\end{cor}
\vspace{2mm}

If $K=\mathbb R$, then $\mathcal H_1=\{-1, 1\}$ from Lemma~\ref{lem:equiv1} (ii). Then
Theorem~\ref{thm:char-not=2} implies the following.

\begin{cor}\label{cor:K=R}% Corollary 7. 2
Suppose $K=\mathbb R$.  Up to isomorphism, 2-dimensional endo-commutative straight algebras of type II$_1$ over $K$ are classified into 4 families $S(0, 1, 1, 0, -1, 2)$, $S(0, 4, -4, -4, 4, 0)$, $\left\{S(0, \varepsilon t, t, 0, \delta t, 0)\right\}_{t, \varepsilon, \delta=\pm1}$ and 
\[
 \{S(0, (1+t)^2/4, (t^2-1)/4, 1, (1-t^2)/4, t)\}_{t\in K^* \backslash \{\pm1\}}.
\]
\end{cor}

 If $K=\mathbb Q$, then $\mathcal H_1$ is countably infinite from Lemma~\ref{lem:equiv1} (iii), and hence the family $\left\{S(0, \varepsilon t, t, 0, \delta t, 0)\right\}_{t\in\mathcal H_1, \varepsilon, \delta=\pm1}$ in Theorem~\ref{thm:char-not=2} is countably infinite.

\begin{thm}\label{thm:char=2}% Theorem 7. 2
Suppose ${\rm{char}}\, K=2$.  Up to isomorphism, 2-dimensional endo-commutative straight algebras of type II$_1$ over $K$ are classified into 4 families
\[
\{S(0, t, t, 0, t, 1)\}_{t\in\mathcal H_2}, \{S(0, t, t, 0, t, 0)\}_{t\in\mathcal H_3}, \{S(0, t, t, t, t, 0)\}_{ t\in \mathcal H_4}.
\]
and 
\[
\left\{S\left(0, \frac{t^2}{1+t^2}, \frac{t}{1+t^2}, 1, \frac{t}{1+t^2}, 1\right)\right\}_{t\in K^*\backslash \{1\}}
\]
with multiplication tables on a linear basis $\{e, f\}$ defined by
\[
\begin{pmatrix}f&te\\te+f&tf\end{pmatrix}, \begin{pmatrix}f&te\\te&tf\end{pmatrix}, \begin{pmatrix}f&t(e+f)\\te&tf\end{pmatrix}\, \, {\rm{and}}\, \, \begin{pmatrix}f&\frac{t}{1+t^2}e+f\\\frac{t}{1+t^2}e+f&\frac{t^2}{1+t^2}f\end{pmatrix},
\]
respectively.
\end{thm}
\vspace{2mm}

If $K$ is quadratically non-extendable and ${\rm{char}}\, K=2$, then $\mathcal H_2=\{1\}$, $\mathcal H_3=\{1\}$ and $\mathcal H_4=\{1\}$ as observed in (II), (III) and (IV) in the previous section, respectively. Then Theorem~\ref{thm:char=2} implies the following.

\begin{cor}\label{cor:quad-non-extendable}% Corollary 7. 3
Suppose that $K$ is quadratically non-extendable and ${\rm{char}}\, K=2$.  Up to isomorphism, 2-dimensional endo-commutative straight algebras of type II$_1$ over $K$ are classified into 4 families $S(0, 1, 1, 0, 1, 1), S(0, 1, 1, 0, 1, 0), S(0, 1, 1, 1, 1, 0)$ and 
$\left\{S\left(0, \frac{t^2}{1+t^2}, \frac{t}{1+t^2}, 1, \frac{t}{1+t^2}, 1\right)\right\}_{t\in K^*\backslash \{1\}}$.
\end{cor}
\vspace{2mm}

If $K=\mathbb F_2(X)$, then both $\mathcal H_2$ and $\mathcal H_4$ are countably infinite from Lemmas~\ref{lem:equiv2} and \ref{lem:equiv5}, hence $\{S(0, t, t, 0, t, 1)\}_{t\in\mathcal H_2}$ and $\{S(0, t, t, t, t, 0)\}_{ t\in \mathcal H_4}$ in Theorem~\ref{thm:char=2} are both countably infinite.  Also since $\mathcal H_3=\{1, X\}$ from Lemma~\ref{lem:equiv3}, we see that  
$\{S(0, t, t, 0, t, 0)\}_{t\in\mathcal H_3}$ in Theorem~\ref{thm:char=2} consists of two algebras $S(0, 1, 1, 0, 1, 0)$ and $S(0, X, X, 0, X, 0)$.  However if $K=\mathbb F_2((X_i)_{i\in\mathbb N})$, then $\mathcal H_3$ is countably infinite from Lemma~\ref{lem:equiv3-infinity}, hence $\{S(0, t, t, 0, t, 0)\}_{t\in\mathcal H_3}$ in Theorem~\ref{thm:char=2} is countably infinite.
\vspace{2mm}

The rest of the paper is devoted to proving our main theorems. The strategy is outlined as follows.
First, we provide classifications of the algebras in $(\mathcal{ECS}_{011})_1$, $(\mathcal{ECS}_{011})_3$, and $(\mathcal{ECS}_{011})_4$, respectively.
Note that there is no need to consider the classification of $(\mathcal{ECS}_{011})_2$, since it is empty by Lemma~\ref{lem:II_1_2}.
Second, we explore isomorphism relationships between the algebras within any pair of the three classes: $(\mathcal{ECS}_{011})_1$, $(\mathcal{ECS}_{011})_3$, and $(\mathcal{ECS}_{011})_4$.

\section{Classification of $(\mathcal{ECS}_{011})_{1}$}\label{sec:ECS1} %section 8

Define
\[
(\mathcal{ECS}_{011})_{1, 1}=\{S(0, a, a, 0, -a, d) :  a, d\ne0\}
\]
and
\[
(\mathcal{ECS}_{011})_{1, 2}=\{S(0, \varepsilon a, a, 0, \delta a, 0) : a\ne0, \varepsilon, \delta=\pm1\}.
\]
By Lemma~\ref{lem:II_1_1}, we have
\[
(\mathcal{ECS}_{011})_{1}=(\mathcal{ECS}_{011})_{1, 1}\sqcup(\mathcal{ECS}_{011})_{1, 2}.
\]

(I) Classification of $(\mathcal{ECS}_{011})_{1, 1}$
\vspace{2mm}

Let $a, d, a', d'\in K^*$ be arbitrary.  By Lemma~\ref{lem:iso}, we see that $S(0, a, a, 0, -a, d) \cong S(0, a', a', 0, -a', d')$ iff there are $x, y, z, w\in K$ with $\left|\begin{matrix}x&y\\z&w\end{matrix}\right|\ne0$ such that (\ref{eq:st-iso}) holds.  In this case, since $p=p'=0, q=a, b=0, c=-a, q'=a', b'=0$ and $c'=-a'$, (\ref{eq:st-iso}) is rewritten as 
\[
\left\{\begin{array}{@{\,}lll}
0=z, \\
x^2+a'y^2+d'xy=w, \\
0=az, \\
z^2+a'w^2+d'zw=aw, \\
a'xw-a'yz=ax, \\
xz+a'yw+d'yz=ay, \\
a'yz-a'xw=-ax+dz, \\
xz+a'yw+d'xw=-ay+dw,
\end{array} \right. 
\]
which is further rewritten as
\[
(\flat)\equiv\left\{\begin{array}{@{\,}lll}
z=0\\
x^2+a'y^2+d'xy=w\\
a'w=a\\
2ay+d'xw=dw
\end{array} \right. 
\]
because $x, w\ne0$. 
\vspace{2mm} 

(I-1) The case where ${\rm{char}}\, K\ne2$.  In this case, by the third and 4th equations of $(\flat)$, we obtain $y=(dw-d'xw)/2a=(d-d'x)/2a'$.  Substituting this $y$ and $w=a/a'$ into the second equation of $(\flat)$, we obtain
\[
x^2+\frac{(d-d'x)^2}{4a'}+\frac{d'x(d-d'x)}{2a'}=\frac{a}{a'},
\]
hence
\[
4a'x^2+d^2-2dd'x+d'^2x^2+2dd'x-2d'^2x^2=4a,
\]
so
\[
4a'x^2+d^2-d'^2x^2=4a.
\]
Thus we obtain $(4a'-d'^2)x^2=4a-d^2$.   Therefore $(\flat)$ is rewritten as  
\[
(\flat')\equiv\left\{\begin{array}{@{\,}lll}
z=0\\
w=a/a'\\
(4a'-d'^2)x^2=4a-d^2\\
2a'y+d'x=d.
\end{array} \right. 
\]
Therefore, we have $S(0, a, a, 0, -a, d) \cong S(0, a', a', 0, -a', d')$ iff there are $x, y, z, w\in K$ with $xw\ne0$ such that ($\flat'$) holds.

Put
 \[
 \mathcal S_0=\{S(0, a, a, 0, -a, d) : a, d\ne0, 4a-d^2=0\}
 \]
and
\[
 \mathcal S_1=\{S(0, a, a, 0, -a, d) : a, d\ne0, 4a-d^2\ne0\}.
 \]
Then we have $(\mathcal{ECS}_{011})_{1, 1}=\mathcal S_0\sqcup\mathcal S_1$.
\vspace{1mm}

(I-1-1) We see that each algebra in $\mathcal S_0$ is not isomorphic to any algebra in $\mathcal S_1$.  Otherwise, $S(0, a, a, 0, -a, d)\cong S(0, a', a', 0, -a', d')$ for some $S(0, a, a, 0, -a, d)\in\mathcal S_0$ and $S(0, a', a', 0, -a', d')\in\mathcal S_1$.  Then there are $x, y, z, w\in K$ with  $x\ne0$ such that $(\flat')$ holds.  By the third equation of $(\flat')$, we have $x^2=(4a-d^2)/(4a'-d'^2)=0/(4a'-d'^2)=0$, so $x=0$, a contradiction.  
\vspace{3mm}

(I-1-2)  We see that any algebra in $\mathcal S_0$ is isomorphic to $S(0, 1, 1, 0, -1, 2)$.  In fact, take $S(0, a, a, 0, -a, d)\in \mathcal S_0$ arbitrarily.  By $(\flat')$, we can easily see that $S(0, a, a, 0, -a, d)\cong S(0, 1, 1, 0, -1, 2)$ with a transformation matrix $\begin{pmatrix}1&(d-2)/2\\0&a\end{pmatrix}$.  

(I-1-3) Define
\[
\mathcal S'_1:=\left\{S\left(0, \frac{t+1}{4}, \frac{t+1}{4}, 0, -\frac{t+1}{4}, 1\right) : t\in K^*\, \backslash\, \{-1\}\right\},
\]
and hence $\mathcal S'_1$ is a subfamily of $\mathcal S_1$.   We assert that any algebra in $\mathcal S_1$ is isomorphic to some algebra in $\mathcal S'_1$.  To show this, define
\[
\mathcal S_1(t):=\{S(0, a, a, 0, -a, d) : a, d\ne0, 4a-d^2=t\}
\]
for each $t\in K^*$, and hence $\mathcal S_1=\bigcup_{t\in K^*}\mathcal S_1(t)$.   Note that $\mathcal S_1(t)\ne\emptyset $ for each $t\in K^*\backslash\{-1\}$ because $4a-d^2=t$ when $t\ne-1, a=(t+1)/4$ and $d=1$.  However, it is possible that $\mathcal S_1(-1)=\emptyset$.  We have the following.
\vspace{2mm}

(a)   For each $t\in K^*$, all algebras in $\mathcal S_1(t)$ are isomorphic to each other.  In fact, if $S(0, a, a, 0, -a, d), S(0, a', a', 0, -a', d')\in\mathcal S_1(t)$, we easily see from $(\flat')$ that these algebras are isomorphic to each other with a transformation matrix $\begin{pmatrix}1&(d-d')/2a'\\0&a/a'\end{pmatrix}$.  
\vspace{2mm}

Now let $S\in \mathcal S_1$ be arbitrary, and hence there is $t\in K^*$ with $S\in\mathcal S_1(t)$ because $\mathcal S_1=\bigcup_{t\in K^*}\mathcal S_1(t)$.
\vspace{2mm}

(b)  The case where $t\in K^*\backslash\{-1\}$.  Put  $S_t=S(0, \frac{t+1}{4}, \frac{t+1}{4}, 0, -\frac{t+1}{4}, 1)$.  Since $S_t\in\mathcal S_1(t)$, it follows from (a) that $S\cong S_t$.  Of course,  $S_t$ is in $\mathcal S'_1$.
\vspace{2mm}

(c)  The case where $t=-1$.  Note that 
\[
\mathcal S_1(-1)=\{S(0, (d^2-1)/4, (d^2-1)/4, 0, -(d^2-1)/4, d) : d\ne0, 1, -1\}, 
\]
hence $\mathcal S_1(-1)=\emptyset$ whenever $\sharp K\le3$.  Suppose $\sharp K\ge4$.   Take $\alpha\in K\backslash\{0, 1, -1\}$ and put $S_{\alpha}=S(0, \frac{\alpha^2-1}{4}, \frac{\alpha^2-1}{4}, 0, -\frac{\alpha^2-1}{4}, \alpha)$. Since $S_{\alpha}\in \mathcal S_1(-1)$, it follows from (a) that $S\cong S_{\alpha}$.  Put $t_0=-\alpha^2$, hence $t_0\in K^*\backslash\{-1\}$.  Moreover, we see easily from $(\flat')$ that $S(0, \frac{t_0+1}{4}, \frac{t_0+1}{4}, 0, -\frac{t_0+1}{4}, 1)\cong S_{\alpha}$ with a transformation matrix $\begin{pmatrix}\alpha&-2\\0&-1\end{pmatrix}$.  Then $S\cong S(0, \frac{t_0+1}{4}, \frac{t_0+1}{4}, 0, -\frac{t_0+1}{4}, 1)\in\mathcal S'_1$.

By (b) and (c), we see that $S$ is isomorphic to some algebra in $\mathcal S'_1$, that is, our assertion holds.
\vspace{3mm}

(I-1-4) Recall that for any $t, t'\in K^*$, we write $t\sim_1t'$ if $t/t'\in(K^*)^2$ (see (I) in \S6).  Let $t, t'\in K^*\, \backslash\, \{-1\}$ and put
\[
S_t=S\left(0, \frac{t+1}{4}, \frac{t+1}{4}, 0, -\frac{t+1}{4}, 1\right)\, \, {\rm{and}}\, \, S_{t'}=S\left(0, \frac{t'+1}{4}, \frac{t'+1}{4}, 0, -\frac{t'+1}{4}, 1\right).
\]
Then we see that $S_t\cong S_{t'}$ iff $t\sim_1 t'$ as follows.  Suppose first $S_t\cong S_{t'}$.  Then there are $x, y, z, w\in K$ with $\left|\begin{matrix}x&y\\z&w\end{matrix}\right|\ne0$ such that 
\[
(\flat')\equiv\left\{\begin{array}{@{\,}lll}
z=0, \\
w=a/a', \\
(4a'-d'^2)x^2=4a-d^2, \\
2a'y+d'x=d.
\end{array} \right. 
\]
holds, where $\left\{\begin{array}{@{\,}lll}a=(t+1)/4, \\d=1\end{array} \right.$ and $\left\{\begin{array}{@{\,}lll}a'=(t'+1)/4, \\d'=1.\end{array} \right.$  By the third equation of $(\flat')$, we have $t'x^2=t$.  But since $\left|\begin{matrix}x&y\\z&w\end{matrix}\right|\ne0$ and $z=0$, it follows that $x\ne0$, and hence $t/t'=x^2\in (K^*)^2$, that is, $t\sim_1t'$.

Conversely, suppose that $t\sim_1t'$.  Then there is $x_0\in K^*$ with $t=t'x_0^2$. Define
\[
a=\frac{t+1}{4},d=1, a'=\frac{t'+1}{4}, d'=1, y_0=\frac{d-d'x_0}{2a'}, z_0=0\, \, {\rm{and}}\, \, w_0=\frac{a}{a'}.
\]
Then $(x, y, z, w):=(x_0, y_0, z_0, w_0)$ satisfy $(\flat')$ because 
\[
\left|\begin{matrix}x_0&y_0\\z_0&w_0\end{matrix}\right|=x_0w_0=\frac{x_0a}{a'}\ne0, (4a'-d'^2)x_0^2=t'x_0^2=t=4a-d^2
\]
and
\[
2a'y_0+d'x_0=d-d'x_0+d'x_0=d.
\]
Therefore we obtain  $S_t\cong S_{t'}$ with a transformation matrix $\begin{pmatrix}x_0&y_0\\z_0&w_0\end{pmatrix}$.

Let $\mathcal H'_1$ be an arbitrary complete representative system of  $K^*\backslash\{-1\}/\sim_1$.   Therefore, by (I-1-1), (I-1-2), (I-1-3) and (I-1-4), we have the following.
\vspace{3mm}

\begin{pro}\label{pro:8.1} % Proposition 8.1 
Suppose ${\rm{char}}\, K\ne2$.  Then all algebras in $(\mathcal{ECS}_{011})_{1, 1}$ are classified into 2 families $S(0, 1, 1, 0, -1, 2)$ and $ \left\{S\left(0, \frac{t+1}{4}, \frac{t+1}{4}, 0, -\frac{t+1}{4}, 1\right)\right\}_{ t\in \mathcal H'_{1}}$ up to isomorphism.
\end{pro}
\vspace{3mm}

(I-2) The case where ${\rm{char}}\, K=2$.  In this case, the 4th equation of $(\flat)$ becomes $d'xw=dw$, and hence $(\flat)$ is rewritten as  
\[
(\flat'')\equiv\left\{\begin{array}{@{\,}lll}
z=0, \\
a'y^2+dy=a/a'+d^2/d'^2, \\
w=a/a', \\
x=d/d'.
\end{array} \right. 
\]
Moreover, we see that $\mathcal S_0=\emptyset$ and  $\mathcal S_1=\{S(0, a, a, 0, a, d) : a, d\ne0\}$, hence $(\mathcal{ECS}_{011})_{1, 1}=\{S(0, a, a, 0, a, d) : a, d\ne0\}$.  For each $t\in K^*$, define
\[
\mathcal S(t):=\{S(0, a, a, 0, a, d) : a, d\ne0, a/d^2=t\},
\]
and hence we have $(\mathcal{ECS}_{011})_{1, 1}=\bigcup_{t\in K^*}\mathcal S(t)$.  Furthermore, we see easily from $(\flat'')$ that for each $t\in K^*$,
\[
S(0, a, a, 0, a, d) \cong S(0, t, t, 0, t, 1)\quad(\forall S(0, a, a, 0, a, d)\in \mathcal S(t))
\]
with a transformation matrix $\begin{pmatrix}d&0\\0&a/t\end{pmatrix}$.  Then we have 
\begin{equation}\label{eq:t}% (9)
\forall S(0, a, a, 0, a, d)\in(\mathcal{ECS}_{011})_{1, 1}, \exists t\in K^* : S(0, t, t, 0, t, 1)\cong S(0, a, a, 0, a, d).
\end{equation}

Recall that for any $t, t'\in K^*$, we write $t\sim_{2}t'$ if $t+t'\in K_{\iota^2+\iota}$ (see (II) in Section~\ref{sec:equiv}).  If $t, t'\in K^*$, then 
\begin{equation}\label{eq:equiv2} % (10)
S(0, t, t, 0, t, 1)\cong S(0, t', t', 0, t', 1)\Leftrightarrow t\sim_{2}t'.
\end{equation}
Indeed, let $t, t'\in K^*$.  Then we see from ($\flat''$) that $S(0, t, t, 0, t, 1)\cong S(0, t', t', 0, t', 1)$ iff there are $x, y, z, w\in K$ with $\left|\begin{matrix}x&y\\z&w\end{matrix}\right|\ne0$ such that 
\[
(\flat''')\equiv\left\{\begin{array}{@{\,}lll}
z=0, \\
t'y^2+y=t/t'+1, \\
w=t/t', \\
x=1
\end{array} \right. 
\]
holds.  The second equation of $(\flat''')$ is rewritten as $(t'y)^2+t'y+t+t'=0$, hence $t\sim_{2}t'$.  Then we obtain $S(0, t, t, 0, t, 1)\cong S(0, t', t', 0, t', 1)\Rightarrow t\sim_{2}t'$.  

Conversely assume $t\sim_{2}t'$.  Then there is $\alpha\in K$ such that $t+t'=\alpha^2+\alpha$.  Put $x=1, y=\alpha/t', z=0$ and $w=t/t'$, hence these scalars satisfy ($\flat'''$) because $t'y^2+y=\alpha^2/t'+\alpha/t'=t/t'+1$.  Then  $S(0, t, t, 0, t, 1)\cong S(0, t', t', 0, t', 1)$ with a transformation matrix $\begin{pmatrix}1&\alpha/t'\\0&t/t'\end{pmatrix}$.  Thus we have shown that (\ref{eq:equiv2}) is true.

Let $\mathcal H_2$ denote a complete representative system of $K^*/\sim_{2}$.  From (\ref{eq:t}) and (\ref{eq:equiv2}), we obtain the following.

\begin{pro}\label{pro:8.2} % Proposition 8. 2
Suppose ${\rm{char}}\, K=2$.  Then all algebras in $(\mathcal{ECS}_{011})_{1, 1}$ are classified into the family $\{S(0, t, t, 0, t, 1)\}_{t\in\mathcal H_2}$ up to isomorphism.
\end{pro}

(II) Classification of $(\mathcal{ECS}_{011})_{1, 2}$
\vspace{2mm}

Let $a, a'\in K^*$ and $\varepsilon, \delta, \varepsilon', \delta'\in \{1, -1\}$ be arbitrary.  By Lemma~\ref{lem:iso}, we see that $S(0, \varepsilon a, a, 0, \delta a, 0)\cong S(0, \varepsilon' a', a', 0, \delta' a', 0)$ iff there are $x, y, z, w\in K$ with $\left|\begin{matrix}x&y\\z&w\end{matrix}\right|\ne0$ such that (\ref{eq:st-iso}) holds with $p=p'=0, q=\varepsilon a, b=0, c=\delta a, q'=\varepsilon'a', b'=0$ and $c'=\delta'a'$ and $d=d'=0$.  In this case, (\ref{eq:st-iso}) is rewritten as 
\[
(\flat)\equiv\left\{\begin{array}{@{\,}lll}
a'(1+\delta')xy=z, \cdots\flat_1\\
x^2+\varepsilon'a'y^2=w, \cdots\flat_2\\
a'(1+\delta')zw=\varepsilon az, \cdots\flat_3\\
z^2+\varepsilon'a'w^2=\varepsilon aw, \cdots\flat_4\\
a'xw+\delta'a'yz=ax, \cdots\flat_5\\
xz+\varepsilon'a'yw=ay, \cdots\flat_6\\
a'yz+\delta'a'xw=\delta ax, \cdots\flat_7\\
xz+\varepsilon'a'yw=\delta ay. \cdots\flat_8
\end{array} \right. 
\]

(II-1) The case where ${\rm{char}}\, K\ne2$.  Suppose $S(0, \varepsilon a, a, 0, \delta a, 0)\cong S(0, \varepsilon' a', a', 0, \delta' a', 0)$.   Then there are $x, y, z, w\in K$ with $\left|\begin{matrix}x&y\\z&w\end{matrix}\right|\ne0$ such that ($\flat$) holds from the above observation.  By $\flat_5-\delta'\times\flat_7$, we have $ax=\delta'\delta ax$, hence $x=\delta'\delta x$ because $a\ne0$.  If $\delta'\ne\delta$, then $\delta'\delta=-1$, hence $x=0$, and so $z=0$ by $\flat_1$. This implies $\left|\begin{matrix}x&y\\z&w\end{matrix}\right|=0$, a contradiction.  Thus we obtain $\delta'=\delta$, and hence $\flat_5\Leftrightarrow\flat_7$.  Moreover, $\flat_6-\flat _8$ yields $y=\delta y$.  Then $(\flat)$ is rewritten as 
\[
(\flat)_1\equiv\left\{\begin{array}{@{\,}lll}
a'(1+\delta)xy=z, \\
x^2+\varepsilon'a'y^2=w, \\
a'(1+\delta)zw=\varepsilon az, \\
z^2+\varepsilon'a'w^2=\varepsilon aw, \\
a'xw+\delta a'yz=ax, \\
xz+\varepsilon'a'yw=ay, \\
y=\delta y.
\end{array} \right. 
\]
(i) The case where $\delta=-1$.  The 7th equation of ($\flat)_1$ becomes $y=0$.  Then we can easily see that $(\flat)_1$ is rewritten as 
\[
\left\{\begin{array}{@{\,}lll}
y=z=0, \\
x^2=w, \\
\varepsilon'=\varepsilon, \\
a'w=a.
\end{array} \right. 
\]
Therefore, we obtain $\varepsilon'=\varepsilon$ and $a\sim_1a'$ because $a/a'=w=x^2\in(K^*)^2$.
\vspace{2mm}

\noindent (ii) The case where $\delta=1$.  In this case, $(\flat)_1$ is rewritten as 
\[
(\flat)_{1,1}\equiv\left\{\begin{array}{@{\,}lll}
2a'xy=z, \\
x^2+\varepsilon'a'y^2=w, \\
2a'zw=\varepsilon az, \\
z^2+\varepsilon'a'w^2=\varepsilon aw, \\
a'xw+a'yz=ax, \\
xz+\varepsilon'a'yw=ay.
\end{array} \right. 
\]
Then we show that $(\flat)_{1,1}$ is rewritten as 
\[
(\flat)_{1,1, 1}\equiv\left\{\begin{array}{@{\,}lll}
y=z=0, \\
x^2=w, \\
\varepsilon'=\varepsilon, \\
a'w=a\end{array} \right. 
\, \, {\rm{or}}\, \, 
(\flat)_{1,1, 2}\equiv\left\{\begin{array}{@{\,}lll}
z=2a'xy\ne0, \\
w=\varepsilon a/2a', \\
2a'x^2+2\varepsilon'a'^2y^2=\varepsilon a, \\
16a'^3x^2y^2+\varepsilon'a^2=2a^2, \\
\varepsilon a+4a'^2y^2=2a, \\
4a'x^2+\varepsilon'\varepsilon a=2a.
\end{array} \right. 
\]
In fact, if $z=0$, then $x\ne0$ and $w\ne0$ because $\left|\begin{matrix}x&y\\z&w\end{matrix}\right|\ne0$, hence $y=0$ by the first equation of $(\flat)_{1,1}$.  Then we can easily see that $(\flat)_{1, 1}$ is rewritten as $(\flat)_{1, 1, 1}$.  Otherwise, we get $w=\varepsilon a/2a'$ from the third equation of $(\flat)_{1, 1}$.  So substituting this $w$ and $z=2a'xy$ into other equations of $(\flat)_{1, 1}$, we can easily see that $(\flat)_{1, 1}$ is rewritten as $(\flat)_{1, 1, 2}$.  However, since  
\[
x^2=\frac{(2-\varepsilon'\varepsilon)a}{4a'}\quad{\rm{and}}\quad y^2=\frac{(2-\varepsilon)a}{4a'^2}
\]
from the 5th and 6th equations of $(\flat)_{1, 1, 2}$, we see that the third and 4th equations of $(\flat)_{1, 1, 2}$ yield $(1+\varepsilon')(1-\varepsilon)=0$ by a simple calculation.  Also, since
\[
\left|\begin{matrix}x&y\\z&w\end{matrix}\right|=\frac{\varepsilon ax}{2a'}-2a'xy^2=\frac{x}{2a'}\left(\varepsilon a- 4a'^2\times\frac{(2-\varepsilon)a}{4a'^2}\right)=\frac{xa}{2a'}\times(2\varepsilon-2)\ne0,
\]
it follows that $\varepsilon\ne1$.  This implies  $\varepsilon=\varepsilon'=-1$ because $(1+\varepsilon')(1-\varepsilon)=0$.  Therefore $(\flat)_{1, 1, 2}$ is rewritten as 
\[
\left\{\begin{array}{@{\,}lll}
z=2a'xy\ne0, \\
w=\varepsilon a/2a', \\
x^2=(2-\varepsilon'\varepsilon)a/4a', \\
y^2=(2-\varepsilon)a/4a'^2, \\
\varepsilon=\varepsilon'=-1,
\end{array} \right. 
\]
that is, 
\[
(\flat)_{1,1, 3}\equiv\left\{\begin{array}{@{\,}lll}
z=2a'xy\ne0, \\
w=-a/2a', \\
x^2=a/4a', \\
y^2=3a/4a'^2, \\
\varepsilon=\varepsilon'=-1.
\end{array} \right. 
\]
In the case where $(\flat)_{1, 1, 1}$ holds, we have $\varepsilon=\varepsilon'$ and $a/a'=w=x^2\in(K^*)^2$, that is, $a\sim_1a'$.  In the case where $(\flat)_{1, 1, 2}$ holds, we have $\varepsilon=\varepsilon'$ and $a/a'=4x^2=(2x)^2\in(K^*)^2$ from $(\flat)_{1,1, 3}$, that is, $a\sim_1a'$.  In any case, we obtain  $\left\{\begin{array}{@{\,}lll}\varepsilon=\varepsilon'\\\delta=\delta'\end{array} \right.$ and $a\sim_1a'$.  

Conversely, suppose that $\left\{\begin{array}{@{\,}lll}\varepsilon=\varepsilon'\\\delta=\delta'\end{array} \right.$ and $a\sim_1a'$.  Then there is $x_0\in K^*$ with $x_0^2=a/a'$.  Define $y_0=0, z_0=0$ and $w_0=a/a'$, hence $\left|\begin{matrix}x_0&y_0\\z_0&w_0\end{matrix}\right|=x_0a/a'\ne0$.  Also it is obvious that $(x, y, z, w):=(x_0, y_0, z_0, w_0)$ satisfies $(\flat)_{1}$.  Therefore $S(0, \varepsilon a, a, 0, \delta a, 0)\cong S(0, \varepsilon' a', a', 0, \delta' a', 0)$ with a transformation matrix $\begin{pmatrix}x_0&y_0\\z_0&w_0\end{pmatrix}$.  
Summarizing the above, we have the following.

\begin{lem}\label{lem:8.1} % Lemma 8.1
Let $a, a'\in K^*$ and $\varepsilon, \delta, \varepsilon', \delta'\in \{1, -1\}$.  Then $S(0, \varepsilon a, a, 0, \delta a, 0)\cong S(0, \varepsilon' a', a', 0, \delta' a', 0)$ iff $\left\{\begin{array}{@{\,}lll}\varepsilon=\varepsilon'\\\delta=\delta'\end{array} \right.$ and $a\sim_1a'$. 
\end{lem}

Recall that $\mathcal H_1$ denote a complete representative system of $K^*/\sim_1$.  Then we have the following.

\begin{pro} \label{pro:8.3}% Proposition 8.3
Suppose ${\rm{char}}\, K\ne2$.  Then all algebras in $(\mathcal{ECS}_{011})_{1, 2}$ are classified into the family $\{S(0, \varepsilon t, t, 0, \delta t, 0)\}_{t\in\mathcal H_1, \varepsilon, \delta=\pm1}$ up to isomorphism.
\end{pro}

\begin{proof}
This immediately follows from Lemma~\ref{lem:8.1} .
\end{proof}

(II-2) The case where ${\rm{char}}\, K=2$.  Recall that for any $t, t'\in K^*$, we write $t\sim_3t'$ if there are $x\in K^*$ and $y\in K$ such that $t'x^2+y^2+t=0$ (see (III) in \S6).  Let $\mathcal H_3$ denote a complete representative system of $K^*/\sim_3$.  Then we have the following.

\begin{pro} \label{pro:8.4}% Proposition 8.4
Suppose ${\rm{char}}\, K=2$.  Then all algebras in $(\mathcal{ECS}_{011})_{1, 2}$ are classified into the family $\{S(0, t, t, 0, t, 0)\}_{t\in\mathcal H_3}$ up to isomorphism.
\end{pro}

\begin{proof}
Since  ${\rm{char}}\, K=2$, it follows that $(\mathcal{ECS}_{011})_{1, 2}=\{S(0, a, a, 0, a, 0) : a\in K^*\}$.  Take $a, a'\in K^*$ arbitrarily.  Then $(\flat)$ is rewritten as 
\[
(\flat')\equiv\left\{\begin{array}{@{\,}lll}
0=z, \\
x^2+a'y^2=w, \\
0=az, \\
z^2+a'w^2=aw, \\
a'xw+a'yz=ax, \\
xz+a'yw=ay, \\
a'yz+a'xw=ax, \\
xz+a'yw=ay.
\end{array} \right. 
\]
But since $\left|\begin{matrix}x&y\\z&w\end{matrix}\right|\ne0$, we see $x, w\ne0$ from the first equation of $(\flat')$.  Then we can easily see that $(\flat')$ is rewritten as $(\flat'')\equiv\left\{\begin{array}{@{\,}lll}
z=0, \\
x^2+a'y^2=w, \\
a'w=a.
\end{array} \right.$  Therefore, we see that 
\begin{equation}\label{eq:equiv3} % (11)
S(0, a, a, 0, a, 0)\cong S(0, a', a', 0, a', 0)\Leftrightarrow a\sim_3 a'.
\end{equation} 
In fact, suppose $S(0, a, a, 0, a, 0)\cong S(0, a', a', 0, a', 0)$. Then there are $x, y, z, w\in K$ with $x\ne0$ such that ($\flat''$) holds.  By the second and third equations of $(\flat'')$, we have $a'x^2+(a'y)^2+a=0$, and hence $a\sim_3a'$.  Conversely, suppose $a\sim_3a'$.  Then, there are $x_0\in K^*$ and $y_0\in K$ such that $a'x_0^2+y_0^2+a=0$, and hence $(x, y, z, w):=(x_0, y_0/a', 0, a/a')$ obviously satisfies ($\flat'')$.  Then $S(0, a, a, 0, a, 0)\cong S(0, a', a', 0, a', 0)$ with a transformation matrix $\begin{pmatrix}x_0&y_0/a'\\0&a/a'\end{pmatrix}$, as claimed.  The desired result easily follows from (\ref{eq:equiv3}).
\end{proof}
 We next investigate the relationship of isomorphism between the algebras in $(\mathcal{ECS}_{011})_{1, 1}$ and $(\mathcal{ECS}_{011})_{1, 2}$
 \vspace{3mm}
 
 (A) The case where ${\rm{char}}\, K\ne2$.
  \vspace{2mm}
 
 (A-1) We assert that $S(0, 1, 1, 0, -1, 2)\ncong S(0, \varepsilon t, t, 0, \delta t, 0)$ for every $t\in K^*$ and $\varepsilon, \delta\in\{1, -1\}$.  In fact, let $t\in K^*$ and $\varepsilon, \delta\in\{1, -1\}$ and suppose $S(0, 1, 1, 0, -1, 2)\cong S(0, \varepsilon t, t, 0, \delta t, 0)$.  Putting
\[
p=0, q=1, a=1, b=0, c=-1, d=2
\]
and
\[
p'=0, q'=\varepsilon t, a'=t, b'=0, c'=\delta t, d'=0
\]
in Lemma~\ref{lem:iso}, we see that  there are $x, y, z, w\in K$ with $\left|\begin{matrix}x&y\\z&w\end{matrix}\right|\ne0$ such that 
\begin{equation}\label{eq:11&12} % (12)
\left\{\begin{array}{@{\,}lll}
(1+\delta )txy=z, \cdots(12-1)\\
x^2+\varepsilon ty^2=w, \cdots(12-2)\\
(1+\delta)tzw=z, \cdots(12-3)\\
z^2+\varepsilon tw^2=w, \cdots(12-4)\\
txw+\delta tyz=x, \cdots(12-5)\\
xz+\varepsilon tyw=y, \cdots(12-6)\\
tyz+\delta txw=-x+2z, \cdots(12-7)\\
xz+\varepsilon tyw=-y+2w\cdots(12-8)
\end{array} \right. 
\end{equation}
holds. By (12-5)$-\delta\times$(12-7),  we obtain $(1+\delta)x=2z$.  Combining this with (12-1), we get $z(2ty-1)=0$, and so either $z=0$ or $2ty=1$.  Assume $z=0$.  Then $x\ne0$ and $w\ne0$ because $\left|\begin{matrix}x&y\\z&w\end{matrix}\right|\ne0$, and so $tw=\varepsilon$ from (12-4).  Combining this with (12-8), we obtain $y=w$.   Then we have from (12-2) that $x^2=w-\varepsilon ty^2=w-\varepsilon tw^2=w-\varepsilon^2w=0$, a contradiction.  Next assume $2ty=1$ and $z\ne0$.  By (12-3), we get $(1+\delta)tw=1$, hence $\delta=1$ and $2tw=1$.  Then we have from (12-1) that $x=z$.  Combining this with (12-5), we obtain $tw+ty=1$, hence $2ty=1$, so $w=y=\frac{1}{2t}$. Thus we get $\left|\begin{matrix}x&y\\z&w\end{matrix}\right|=\left|\begin{matrix}x&y\\x&y\end{matrix}\right|=0$, a contradiction.  Therefore, we have shown that our assertion is true.
  \vspace{3mm}
  
(A-2)  We see that $S\left(0, \frac{t+1}{4}, \frac{t+1}{4}, 0, -\frac{t+1}{4}, 1\right)\cong S(0, t, t, 0, -t, 0)$ holds for all $t\in K^*$.  In fact, let $t\in K^*$, and hence $\begin{pmatrix}1/2&1/2t\\0&(1+t)/4t\end{pmatrix}$ serves as a transformation matrix for this isomorphism.
  \vspace{3mm}

Let $\mathcal H_1$ be a complete representative system of $K^*/\sim_1$.  Then we have the following.

\begin{pro}\label{pro:8.5} % Proposition 8.5
Suppose ${\rm{char}}\, K\ne2$.  Then all algebras in $(\mathcal{ECS}_{011})_{1}$ are classified into 2 families $S(0, 1, 1, 0, -1, 2)$ and $\{S(0, \varepsilon t, t, 0, \delta t, 0)\}_{t\in\mathcal H_1, \varepsilon, \delta=\pm1}$ up to isomorphism
\end{pro}  
\vspace{1mm}

\begin{proof}
   Define
\[
\mathcal H'_1:=\left\{\begin{array}{@{\,}lll}
\mathcal H_1\quad\quad\quad\quad\quad\quad   {\rm{if}}\, \, -1\notin\mathcal H_1, \\
\mathcal H_1\backslash\{-1\}\quad\quad\quad\, \, {\rm{if}}\, \, \widetilde{-1}=\{-1\}, \\
\mathcal H_1\backslash\{-1\}\cup\{\alpha\}\, \,\, {\rm{otherwise}}, 
\end{array} \right. 
\] 
where $\alpha\in\widetilde{-1}\cap K^*$ with $\alpha\ne-1$.  Here $\widetilde{-1}$ denotes the coset of $-1$ with respect to $\sim_1$.  Then we see easily that $\mathcal H'_1$ is a complete representative system of $K^*\backslash\{-1\}/\sim_1$.  Then by Proposition~\ref{pro:8.1}, all algebras in $(\mathcal{ECS}_{011})_{1, 1}$ are classified into 2 families $S(0, 1, 1, 0, -1, 2)$ and $ \left\{S\left(0, \frac{t+1}{4}, \frac{t+1}{4}, 0, -\frac{t+1}{4}, 1\right)\right\}_{ t\in \mathcal H'_{1}}$ up to isomorphism.   By Proposition~\ref{pro:8.3}, all algebras in $(\mathcal{ECS}_{011})_{1, 2}$ are classified into the family $\{S(0, \varepsilon t, t, 0, \delta t, 0)\}_{t\in\mathcal H_1, \varepsilon, \delta=\pm1}$ up to isomorphism.  Then by (A-1), we see that no algebras in $(\mathcal{ECS}_{011})_{1, 2}$ is isomorphic to $S(0, 1, 1, 0, -1, 2)$.  Let $t\in \mathcal H'_1$ be arbitrary.  We assert that $S(0, \frac{t+1}{4}, \frac{t+1}{4}, 0, -\frac{t+1}{4}, 1)$ is isomorphic to some algebra in $(\mathcal{ECS}_{011})_{1, 2}$.  In fact, if either $-1\notin\mathcal H_1$ or $ \widetilde{-1}=\{-1\}$, then $\mathcal  H'_1\subseteq\mathcal H_1$, and hence our assertion is true from (A-2).  Otherwise, we have $-1\in\mathcal H_1$ and $\mathcal H'_1=\mathcal H_1\backslash\{-1\}\cup\{\alpha\}$.  In the case where $t\in\mathcal H_1\backslash\{-1\}$, our assertion is true from (A-2).  In the case where $t=\alpha$, taking $\varepsilon=\varepsilon'=1, \delta=\delta'=-1, a=-1$ and $a'=\alpha$ in Lemma 8.1, we see $S(0, -1, -1, 0, 1, 0)\cong S(0, \alpha, \alpha, 0, -\alpha, 0)$.  Also, $S(0, \frac{\alpha+1}{4}, \frac{\alpha+1}{4}, 0, -\frac{\alpha+1}{4}, 1)\cong S(0, \alpha, \alpha, 0, -\alpha, 0)$ from (A-2).  Then $S(0, \frac{\alpha+1}{4}, \frac{\alpha+1}{4}, 0, -\frac{\alpha+1}{4}, 1)\cong S(0, -1, -1, 0, 1, 0)$ and $S(0, -1, -1, 0, 1, 0)\in(\mathcal{ECS}_{011})_{1, 2}$ because $-1\in\mathcal H_1$, hence our assertion is true.   

Recall that $(\mathcal{ECS}_{011})_{1}=(\mathcal{ECS}_{011})_{1, 1}\sqcup(\mathcal{ECS}_{011})_{1, 2}$.  Therefore, the above observations imply the desired result.
\end{proof}
 
(B) The case where ${\rm{char}}\, K=2$.
\vspace{2mm}

We assert that $S(0, t, t, 0, t, 1)\ncong S(0, t', t', 0, t', 0)$ for all $t, t'\in K^*$.  To show this, suppose $S(0, t, t, 0, t, 1)\cong S(0, t', t', 0, t', 0)$ for some $t, t'\in K^*$.  Putting
\[
p=0, q=t, a=t, b=0, c=t, d=1
\]
and
\[
p'=0, q'=t', a'=t', b'=0, c'=t', d'=0,
\]
we see from Lemma~\ref{lem:iso} that there are $x, y, z, w\in K$ with $\left|\begin{matrix}x&y\\z&w\end{matrix}\right|\ne0$ such that 
\begin{equation}\label{eq:1&0} % (13)
\left\{\begin{array}{@{\,}lll}
2t'xy=z, \cdots(13-1)\\
x^2+t'y^2=w, \cdots(13-2)\\
2t'zw=tz, \cdots(13-3)\\
z^2+t'w^2=tw, \cdots(13-4)\\
t'xw+t'yz=tx, \cdots(13-5)\\
xz+t'yw=ty, \cdots(13-6)\\
t'yz+t'xw=tx+z, \cdots(13-7)\\
xz+t'yw=ty+w\cdots(13-8)
\end{array} \right. 
\end{equation}
holds.  By (13-1), we have $z=0$, hence $w\ne0$ because $\left|\begin{matrix}x&y\\z&w\end{matrix}\right|\ne0$.  Then by (13-4), we have $t'w=t$, hence $w=0$ by (13-8).  This is a contradiction.  Thus our assertion holds.
\vspace{2mm}

Combining (B), Propositions~\ref{pro:8.2} and \ref{pro:8.4}, we have the following.

\begin{pro}\label{pro:8.6} % Proposition 8.6
Suppose ${\rm{char}}\, K=2$.   Then all algebras in $(\mathcal{ECS}_{011})_{1}$ are classified into 2 families $\{S(0, t, t, 0, t, 1)\}_{t\in\mathcal H_2}$ and $\{S(0, t, t, 0, t, 0)\}_{t\in\mathcal H_3}$ up to isomorphism.
\end{pro}  

\section{Classification of $(\mathcal{ECS}_{011})_{3}$}\label{sec:ECS3} %section 9

Recall that $(\mathcal{ECS}_{011})_{3}=\{S(0, -a, a, b, -a, 0) : a, b\ne0\}$ from Lemma~\ref{lem:II_1_3}.  Let $a, b, a', b'\in K^* $ be arbitrary.  Putting
\[
p=0, q=-a, c=-a, d=0, p'=0, q'=-a', c'=-a' \, \, {\rm{and}}\, \, d'=0,
\]
we see from Lemma~\ref{lem:iso} that $S(0, -a, a, b, -a, 0)\cong S(0, -a', a', b', -a', 0)$ iff there are $x, y, z, w\in K$ with $\left|\begin{matrix}x&y\\z&w\end{matrix}\right|\ne0$ such that 
\begin{equation}\label{eq:3-iso} % (14)
\left\{\begin{array}{@{\,}lll}
z=0, \\
x^2-a'y^2+b'xy=w, \\
-a'w^2=-aw, \\
a'xw=ax, \\
-a'yw+b'xw=ay+bw, \\
-a'xw=-ax, \\
-a'yw=-ay
\end{array} \right. 
\end{equation}
holds.  By the first equation of (\ref{eq:3-iso}), we obtain $w\ne0$.  Then $w=a/a'$ by the third equation of (\ref{eq:3-iso}).  Also the 5th equation of (\ref{eq:3-iso}) becomes $b'xw=2ay+bw$.  Therefore (\ref{eq:3-iso}) is rewritten as 
\begin{equation}\label{eq:3-iso-a} % (15)
\left\{\begin{array}{@{\,}lll}
z=0, \\
x^2-a'y^2+b'xy=w, \\
w=a/a', \\
b'xw=2ay+bw.
\end{array} \right. 
\end{equation}
For any $a, b\in K^*$, put $t=b^2/a$, and then we see from (\ref{eq:3-iso-a}) that $S(0, -a, a, b, -a, 0) \cong S(0, -t, t, t, -t, 0)$ with a transformation matrix $\begin{pmatrix}b/t&0\\0&a/t\end{pmatrix}$.  Then any algebra in $(\mathcal{ECS}_{011})_3$ must be isomorphic to $S(0, -t, t, t, -t, 0)$ for some $t\in K^*$.

Let $t, t'\in K^*$.  By (\ref{eq:3-iso-a}), we see that $S(0, -t, t, t, -t, 0)\cong S(0, -t', t', t', -t', 0)$ iff there are $x, y, z, w\in K$ with $\left|\begin{matrix}x&y\\z&w\end{matrix}\right|\ne0$ such that 
\begin{equation}\label{eq:3-iso-b} % (16)
\left\{\begin{array}{@{\,}lll}
z=0, \cdots(16-1)\\
x^2-t'y^2+t'xy=w, \cdots(16-2)\\
w=t/t', \cdots(16-3)\\
x=2y+w\cdots(16-4)
\end{array} \right. 
\end{equation}
holds.  By (16-1), $z=0$ and hence $\left|\begin{matrix}x&y\\z&w\end{matrix}\right|=xw\ne0$, so $x\ne0$.  Substituting (16-3) and (16-4) into (16-2), (16-2) is rewritten as
\begin{equation}\label{eq:16-2-a} % (17)
t'^2(4+t')y^2+tt'(4+t')y+t^2-tt'=0
\end{equation}
from a simple calculation.  Then we have
\begin{equation}\label{eq:non-iso} % (18)
S(0, -t, t, t, -t, 0)\ncong S(0, 4, -4, -4, 4, 0)\, \, (\forall t\in K^*\, \backslash\, \{-4\}).
\end{equation}
In fact, if ${\rm{char}}\, K=2$, then (\ref{eq:non-iso}) obviously holds.  Then we consider the case where ${\rm{char}}\, K\ne2$.  Assume $t\in K^*\, \backslash\, \{-4\}$ and $S(0, -t, t, t, -t, 0)\cong S(0, 4, -4, -4, 4, 0)$.  Replacing $t'$ with $-4$ in (\ref{eq:16-2-a}), we have $t^2-t\times(-4)=0$, that is, $t=-4$, a contradiction.

In the case where ${\rm{char}}\, K\ne2$, recall that for $t, t'\in K^*\backslash\{-4\}$, we write $t\sim_5t'$ if $\frac{t'(4+t)}{t(4+t')}\in(K^*)^2$ (see (V) in \S6).  Let $\mathcal H'_5$ denote a complete representative system of $(K^*\, \backslash\, \{-4\})/\sim_5$.

In the case where ${\rm{char}}\, K=2$, recall that for $t, t'\in K^*$, we write $t\sim_4t'$ if $1/t+1/t'\in K_{\iota^2+\iota}$ (see (IV) in \S6).  Let $\mathcal H_4$ denote a complete representative system of $K^*/\sim_4$.

Then we have the following.
\vspace{2mm}

\begin{pro}\label{pro:9.1} % Proposition 9.1
{\rm{(i)}}  Suppose $ {\rm{char}}\, K\ne2$.  Then all algebras in $(\mathcal{ECS}_{011})_3$ are classified into the family $\{S(0, -t, t, t, -t, 0)\}_{ t\in\{-4\}\cup \mathcal H'_5}$  up to isomorphism.
\vspace{2mm}

{\rm{(ii)}} Suppose $ {\rm{char}}\, K=2$.   Then all algebras in $(\mathcal{ECS}_{011})_3$ are classified into the family $\{S(0, t, t, t, t, 0)\}_{ t\in \mathcal H_4}$ up to isomorphism.
\end{pro}

\begin{proof}
Recall that any algebra in $(\mathcal{ECS}_{011})_3$ is isomorphic to $S(0, -t, t, t, -t, 0)$ for some $t\in K^*$.

(i)  Suppose $ {\rm{char}}\, K\ne2$.  With the help of (\ref{eq:non-iso}), we have only to show that 
\[
S(0, -t, t, t, -t, 0)\cong S(0, -t', t', t', -t', 0)\Leftrightarrow t\sim_5t'
\] 
for each $t, t'\in K^*\, \backslash\, \{-4\}$.   

Let $t, t'\in K^*\, \backslash\, \{-4\}$ and suppose $S(0, -t, t, t, -t, 0)\cong S(0, -t', t', t', -t', 0)$.  Then there is $y\in K$ satisfying (\ref{eq:16-2-a}).  Define $y_0=(t'/t)y$.  Substituting $y=(t/t')y_0$ into (\ref{eq:16-2-a}), we obtain
\[
t'^2(4+t')\times\frac{t^2}{t'^2}y_0^2+tt'(4+t')\times\frac{t}{t'}y_0+t^2-tt'=0,
\]
that is, 
\[
t(4+t')(y_0^2+y_0)+t-t'=0.
\]
Then $y_0\ne-1/2$ (otherwise, $t(4+t')(\frac{1}{4}-\frac{1}{2})+t-t'=0$, hence $t=-4$, a contradiction), and hence
\[
\frac{t'(4+t)}{t(4+t')}=1+4\times\frac{t'-t}{t(4+t')}=1+4\times(y_0^2+y_0)=(2y_0+1)^2\in(K^*)^2,
\]
so $t\sim_5t'$.

Conversely, suppose $t\sim_5t'$.  Then there is $\alpha\in K^*$ such that $\frac{t'(4+t)}{t(4+t')}=\alpha^2$. Define
\[
x=\frac{t\alpha}{t'}, y=\frac{t(\alpha-1)}{2t'}, z=0\, \, {\rm{and}}\, \, w=\frac{t}{t'}.
\]
Then $\left|\begin{matrix}x&y\\z&w\end{matrix}\right|=\frac{t^2\alpha}{t'^2}\ne0$.  It is obvious that scalars $x, y, z, w$ satisfy (16-1) and (16-3).  Also since
\[
2y+w=\frac{2t(\alpha-1)}{2t'}+\frac{t}{t'}=\frac{t\alpha-t+t}{t'}=x,
\]
it follows that scalars $x, y, z, w$ satisfy (16-4).  Moreover since
\begin{align*}
t'^2(4+t')y^2+tt'(4+t')y&=\frac{(4+t')t^2}{4}\alpha^2+\frac{(4+t')t^2}{4}-\frac{t^2(4+t')}{2}\\
&=\frac{tt'(4+t)}{4}-\frac{4t^2+t^2t'}{4}\, \, \left(\because \alpha^2=\frac{t'(4+t)}{t(4+t')}\right)\\
&=tt'-t^2,
\end{align*}
these scalars satisfy (\ref{eq:16-2-a}), that is, (16-2).  Therefore we obtain $S(0, -t, t, t, -t, 0)\cong S(0, -t', t', t', -t', 0)$ with a transformation matrix $\begin{pmatrix}t\alpha/t'&t(\alpha-1)/(2t')\\0&t/t'\end{pmatrix}$.
\vspace{2mm}

(ii)  Suppose $ {\rm{char}}\, K=2$. We have only to show that 
\[
S(0, t, t, t, t, 0)\cong S(0, t', t', t', t', 0)\Leftrightarrow t\sim_4t' 
\]
for each $t, t'\in K^*$.  Let $t, t'\in K^*$ and suppose $S(0, t, t, t, t, 0)\cong S(0, t', t', t', t', 0)$. Then there is $y\in K$ satisfying (\ref{eq:16-2-a}), which is rewritten as $t'^3y^2+tt'^2y+t^2+tt'=0$.  Define $y_0=(t'/t)y$.  Then (\ref{eq:16-2-a}) is rewritten as $y_0^2+y_0+\frac{1}{t'}+\frac{1}{t}=0$.  Then $\frac{1}{t}+\frac{1}{t'}\in K_{\iota^2+\iota}$, that is, $t\sim_4t'$.

Conversely, suppose $t\sim_4t'$, hence there is $\alpha\in K$ such that $\alpha^2+\alpha=\frac{1}{t}+\frac{1}{t'}$.   Define $x=\frac{t}{t'}, y=\frac{t}{t'}\alpha, z=0$ and $w=\frac{t}{t'}$.  Then $\left|\begin{matrix}x&y\\z&w\end{matrix}\right|=\frac{t^2}{t'^2}\ne0$.  It is obvious that scalars $x, y, z, w$ satisfy (16-1), (16-3)   and (16-4).  Also since
\begin{align*}
x^2+t'y^2+t'xy&=\frac{t^2}{t'^2}+\frac{t^2}{t'}(\alpha^2+\alpha)\\
&=\frac{t^2}{t'^2}+\frac{t^2}{t'}\left(\frac{1}{t}+\frac{1}{t'}\right)\, \, \left(\because \alpha^2+\alpha=\frac{1}{t}+\frac{1}{t'}\right)\\
&=2\times\frac{t^2}{t'^2}+\frac{t}{t'}=w,
\end{align*}
these scalars satisfy (16-2).  Therefore we obtain $S(0, t, t, t, t, 0)\cong S(0, t', t', t', t', 0)$ with a transformation matrix $\begin{pmatrix}t/t'&t\alpha/t'\\0&t/t'\end{pmatrix}$.
\end{proof}

\section{Classification of $(\mathcal{ECS}_{011})_{4}$}\label{sec:ECS4} %section 10

(I) The case where ${\rm{char}}K\ne2$.
\vspace{2mm}

Recall that 
\[ 
(\mathcal{ECS}_{011})_4=\{S(0, (b+d)^2/4, (d^2-b^2)/4, b, (b^2-d^2)/4, d) : b, d\ne0, b\ne\pm d\}
 \]
from Lemma~\ref{lem:II_1_4} (i).  Let $b, d, b', d'\in K^* $ be such that $b\ne\pm d$ and $b'\ne\pm d'$.  Putting
\[
p=0, q=(b+d)^2/4, a=(d^2-b^2)/4, c=(b^2-d^2)/4
\]
and
\[
p'=0, q'=(b'+d')^2/4, a'=(d'^2-b'^2)/4, c'=(b'^2-d'^2)/4,
\]
we see from Lemma~\ref{lem:iso} that
\[
S(0, (b+d)^2/4, (d^2-b^2)/4, b, (b^2-d^2)/4, d)\cong S(0, (b'+d')^2/4, (d'^2-b'^2)/4, b', (b'^2-d'^2)/4, d')
\] 
iff there are $x, y, z, w\in K$ with $\left|\begin{matrix}x&y\\z&w\end{matrix}\right|\ne0$ such that 
\begin{equation}\label{eq:4-iso}  % (19)
\left\{\begin{array}{@{\,}lll}
0=z, \\
x^2+(b'+d')^2y^2/4+(b'+d')xy=w, \\
0=0, \\
(b'+d')^2w=(b+d)^2, \\
(d'^2-b'^2)w=d^2-b^2, \\
(b'+d')^2yw/4+b'xw=(d^2-b^2)y/4+bw, \\
(b'^2-d'^2)w=b^2-d^2, \\
(b'+d')^2yw/4+d'xw=(b^2-d^2)y/4+dw
\end{array} \right. 
\end{equation}
holds.   By the 4th and 5th equations of (\ref{eq:4-iso}), we see that 
\[
\frac{(b+d)^2}{(b'+d')^2}=\frac{d^2-b^2}{d'^2-b'^2},
\]
that is, $d'/d=b'/b$.  Put $\tau=d'/d=b'/b$, and hence $\tau\ne0$.  Then (\ref{eq:4-iso}) is rewritten as 
\begin{equation}\label{eq:4-iso-a} % (20)
\left\{\begin{array}{@{\,}lll}
z=0, \\
x^2+(b+d)^2\tau^2y^2/4+(b+d)\tau xy=w, \\
\tau^2w=1, \\
(b+d)^2y/4+b\tau xw=(d^2-b^2)y/4+bw, \\
(b+d)^2y/4+d\tau xw=(b^2-d^2)y/4+dw, \\
\tau=d'/d=b'/b.
\end{array} \right. 
\end{equation}
By the third and 4th equations of (\ref{eq:4-iso-a}), we see that the 5th equation of (\ref{eq:4-iso-a}) becomes $x=-\tau(b+d)y/2+1/\tau$, hence (\ref{eq:4-iso-a}) is further rewritten as
\begin{equation}\label{eq:4-iso-b} % (21)
\left\{\begin{array}{@{\,}lll}
z=0, \\
x^2+(b+d)^2\tau^2y^2/4+(b+d)\tau xy=w, \\
\tau^2w=1, \\
(b+d)^2y/4+b\tau xw=(d^2-b^2)y/4+bw, \\
x=-\tau(b+d)y/2+1/\tau, \\
\tau=d'/d=b'/b.
\end{array} \right. 
\end{equation}
For each $t\in K^*\backslash \{\pm1\}$, define 
\[
\mathcal S(t)=\{S(0, (b+d)^2/4, (d^2-b^2)/4, b, (b^2-d^2)/4, d) : b, d\in K^*, d/b=t\}.
\]
Then we have $(\mathcal{ECS}_{011})_4=\bigcup_{t\in K^*\backslash \{\pm1\}}\mathcal S(t)$.  Moreover, we see easily from (\ref{eq:4-iso-b}) that if $t\in K^*\backslash\{\pm1\}$ and $S(0, (b+d)^2/4, (d^2-b^2)/4, b, (b^2-d^2)/4, d)\in\mathcal S(t)$, then
\[
S(0, (b+d)^2/4, (d^2-b^2)/4, b, (b^2-d^2)/4, d)\cong S(0, (1+t)^2/4, (t^2-1)/4, 1, (1-t^2)/4, t)
\]
with a transformation matrix $\begin{pmatrix}b&0\\0&b^2\end{pmatrix}$.  Then any algebra in $(\mathcal{ECS}_{011})_4$ must be isomorphic to $S(0, (1+t)^2/4, (t^2-1)/4, 1, (1-t^2)/4, t)$ for some $t\in K^*\backslash \{\pm1\}$.  Therefore we have the following.
\vspace{2mm}

\begin{pro}\label{pro:10.1} % Proposition 10.1
Suppose ${\rm{char}}K\ne2$.  Then all algebras in $(\mathcal{ECS}_{011})_4$ are classified into the family $\{S(0, (1+t)^2/4, (t^2-1)/4, 1, (1-t^2)/4, t)\}_{t\in K^* \backslash \{\pm1\}}$ up to isomorphism. 
\end{pro}

\begin{proof}
For any $t\in K^*\backslash \{\pm1\}$, define 
\[
S_t:=S(0, (1+t)^2/4, (t^2-1)/4, 1, (1-t^2)/4, t).
\] 
As observed above, any algebra in $(\mathcal{ECS}_{011})_4$ is isomorphic to $S_t$ for some $t\in K^*\backslash \{\pm1\}$.  Also, it is obvious that $\{S_t : t\in K^*\backslash \{\pm1\}\}$ is a subfamily of $(\mathcal{ECS}_{011})_4$. Then we have only to show that if $t, t'\in K^*\backslash \{\pm1\}$ and $S_t\cong S_{t'}$, then $t=t'$.  In fact, if $S_t\cong S_{t'} \, (t, t'\ne0, \pm1)$, then there are $x, y, z, w\in K$ such that (\ref{eq:4-iso-b}) holds with $b=b'=1, d=t$ and $d'=t'$. Therefore we have $t=t'$ by the last equation of (\ref{eq:4-iso-b}).
\end{proof}

(II) The case where ${\rm{char}}K=2$.
\vspace{2mm}

Recall that 
\[
(\mathcal{ECS}_{011})_{4}=\{S(0, q, a, b, a, b) : a, b, q\ne0, q^2+a^2+qb^2=0\}.
\]
from Lemma~\ref{lem:II_1_4} (ii).  Let $q, a, b, q', a', b'\in K^* $ be such that $q^2+a^2+qb^2=0$ and $q'^2+a'^2+q'b'^2=0$.  Putting $p=0, c=a, d=b, p'=0, c'=a'$ and $d'=b'$, we see from Lemma~\ref{lem:iso} that $S(0, q, a, b, a, b)\cong S(0, q', a', b', a', b')$ iff there are $x, y, z, w\in K$ with $\left|\begin{matrix}x&y\\z&w\end{matrix}\right|\ne0$ such that \begin{equation}\label{eq:4-iso-char2} % (22)
\left\{\begin{array}{@{\,}lll}
z=0, \\
x^2+q'y^2=w, \\
q'w^2=qw, \\
a'xw=ax, \\
q'yw+b'xw=ay+bw, \\
a'xw=ax, \\
q'yw+b'xw=ay+bw
\end{array} \right. 
\end{equation}
holds.  By the first equation of (\ref{eq:4-iso-char2}), we obtain $x\ne0$ and $w\ne0$.  Then (\ref{eq:4-iso-char2}) is rewritten as
\begin{equation}\label{eq:4-iso-char2-a} % (23)
\left\{\begin{array}{@{\,}lll}
z=0, \\
x^2+q'y^2=w, \\
q'w=q, \\
a'w=a, \\
qy+b'xw=ay+bw.
\end{array} \right. 
\end{equation}
By the third and 4th equations of (\ref{eq:4-iso-char2-a}), we see that $q/q'=a/a'$.  Put $\tau=q'/q=a'/a$, and hence $\tau\ne0$.  Then (\ref{eq:4-iso-char2-a}) is further rewritten as 
\begin{equation}\label{eq:4-iso-char2-b} % (24)
\left\{\begin{array}{@{\,}lll}
z=0, \\
x^2+\tau qy^2=w, \\
\tau w=1, \\
qy+b'xw=ay+bw, \\
\tau=q'/q=a'/a.
\end{array} \right. 
\end{equation}
Since $q^2+a^2+qb^2=0$ and $\tau^2q^2+\tau^2a^2+\tau qb'^2=0$, it follows that $\tau^2qb^2=\tau qb'^2$, hence $\tau b^2=b'^2$.  By the 4th equation of (\ref{eq:4-iso-char2-b}), we have $x=\{(a+q)y+bw)\}/(b'w)$, and hence
\begin{align*}
x^2+\tau qy^2&=\frac{(a^2+q^2)y^2+b^2w^2+\tau qy^2b'^2w^2}{b'^2w^2}\\
&=\frac{qb^2y^2+b^2w^2+\tau qy^2b'^2w^2}{b'^2w^2}\, \,(\because q^2+a^2=qb^2)\\
&=\frac{q(b^2+b'^2w)y^2+b^2w^2}{b'^2w^2}\, \,(\because \tau w=1)\\
&=\frac{q(b^2+\tau b^2w)y^2+b^2w^2}{ \tau b^2w^2}\, \,(\because \tau b^2=b'^2)\\
&=\frac{q(b^2+b^2)y^2+b^2w^2}{ \tau b^2w^2}\, \,(\because \tau w=1)\\
&=\frac{0+b^2w^2}{ \tau b^2w^2}=\frac{1}{\tau}=w.
\end{align*}
Therefore, the 4th equation of (\ref{eq:4-iso-char2-b}) implies the second equation of (\ref{eq:4-iso-char2-b}).  Then (\ref{eq:4-iso-char2-b}) is rewritten as
\begin{equation}\label{eq:4-iso-char2-c} % (25)
\left\{\begin{array}{@{\,}lll}
z=0, \\
w=q/q'=a/a', \\
qy+b'xw=ay+bw.
\end{array} \right. 
\end{equation}
Substituting the second equation of (\ref{eq:4-iso-char2-c}) into  the third equation of (\ref{eq:4-iso-char2-c}), we obtain $q'y+b'x=a'y+b$.  Then (\ref{eq:4-iso-char2-c}) is
finally rewritten as
\begin{equation}\label{eq:4-iso-char2-d} % (26)
\left\{\begin{array}{@{\,}lll}
z=0, \\
w=q/q'=a/a', \\
b'x+(q'+a')y+b=0.
\end{array} \right. 
\end{equation}
Consequently, we see that if $S(0, q, a, b, a, b), S(0, q', a', b', a', b')\in(\mathcal{ECS}_{011})_4$, then
\[
S(0, q, a, b, a, b)\cong S(0, q', a', b', a', b')\Leftrightarrow\exists x, y, z, w\in K : x\ne0\, \, {\rm{and}}\, \, (\ref{eq:4-iso-char2-d})\, \, {\rm{holds.}}
\]
Let $t\in K^*\, \backslash\, \{1\}$ and define 
\[
\mathcal S(t):=\{S(0, q, a, b, a, b) : q, a, b\in K^*, q^2+a^2+qb^2=0, q/a=t\}.
\]
Then we have $(\mathcal{ECS}_{011})_4=\bigcup_{t\in K^* \backslash\{1\}}\mathcal S(t)$.  In fact, take $S(0, q, a, b, a, b)\in(\mathcal{ECS}_{011})_4$ arbitrarily and put $t=q/a$.  Then $t\in K^*\backslash\{1\}$, otherwise, $qb^2=q^2+a^2=2q^2=0$, a contradiction.  Therefore $S(0, q, a, b, a, b)\in\bigcup_{t\in K^* \backslash\{1\}}\mathcal S(t)$, hence $(\mathcal{ECS}_{011})_4\subseteq\bigcup_{t\in K^* \backslash\{1\}}\mathcal S(t)$.  The opposite inclusion will be obvious.  Moreover, define
\[
S_t:=S\left(0, \frac{t^2}{1+t^2}, \frac{t}{1+t^2}, 1, \frac{t}{1+t^2}, 1\right),
\]
and hence $S_t\in\mathcal S(t)$.  In fact, putting $q=t^2/(1+t^2), a=t/(1+t^2)$ and $b=1$, we see that $q, a, b\in K^*, q/a=t$ and
\[
q^2+a^2+qb^2=\frac{t^4}{1+t^4}+\frac{t^2}{1+t^4}+\frac{t^2}{1+t^2}=\frac{t^4+t^2+t^2(1+t^2)}{1+t^4}=0,
\]
that is, $S_t$ must be in $\mathcal S(t)$.  We assert that $S(0, q, a, b, a, b)\cong S_t$ for all $S(0, q, a, b, a, b)\in\mathcal S(t)$.  To show this, let $S(0, q, a, b, a, b)\in\mathcal S(t)$ be arbitrary and put
\[
x=b, y=0, z=0, w=\frac{a(1+t^2)}{t}, q'=\frac{t^2}{1+t^2}, a'=\frac{t}{1+t^2}\, \, {\rm{and}}\, \, b'=1.
\]
Then it is obvious that $x\ne0, z=0, w=a/a'$ and $b'x+(q'+a')y+b=0$ hold.  Also 
\[
q/q'=\frac{q(1+t^2)}{t^2}=\frac{at(1+t^2)}{t^2}=\frac{a(1+t^2)}{t}=a/a'=w
\]
holds because $q/a=t$.  Therefore, $x, y, z$ and $w$ satisfy (\ref{eq:4-iso-char2-d}), hence we obtain $S(0, q, a, b, a, b)\cong S_t$, so our assertion has been shown.  Since $(\mathcal{ECS}_{011})_4=\bigcup_{t\in K^* \backslash\{1\}}\mathcal S(t)$, we see from our assertion that any algebra in $(\mathcal{ECS}_{011})_4$ is isomorphic to $S_t$ for some $t\in K^*\backslash\{1\}$.  Now let $t, t'\in K^*\backslash\{1\}$ and suppose $S_t\cong S_{t'}$.  Then there are $x, y, z, w\in K$ such that (\ref{eq:4-iso-char2-d}) holds with 
\[
q=\frac{t^2}{1+t^2}, a=\frac{t}{1+t^2}, b=1, q'=\frac{t'^2}{1+t'^2}, a'=\frac{t'}{1+t'^2}\quad{and}\quad b'=1.
\]
By the second equation of (\ref{eq:4-iso-char2-d}), we see easily that $t=t'$.  Therefore we have the following.

\begin{pro}\label{pro:10.2} % Proposition 10.2
Suppose ${\rm{char}}K=2$.  Then all algebras in $(\mathcal{ECS}_{011})_4$ are classified into the family $\left\{S\left(0, \frac{t^2}{1+t^2}, \frac{t}{1+t^2}, 1, \frac{t}{1+t^2}, 1\right)\right\}_{t\in K^*\backslash \{1\}}$ up to isomorphism. 
\end{pro}

\section{Proof of Theorem~\ref{thm:char-not=2}}\label{sec:proof1} % section 11

In this section, we assume that ${\rm{char}}\, K\ne2$.  
\vspace{2mm}

(I) Relationship between $(\mathcal{ECS}_{011})_1$ and $(\mathcal{ECS}_{011})_3$
\vspace{2mm}

Recall that all algebras in $(\mathcal{ECS}_{011})_{1}$ are classified into 2 families $S(0, 1, 1, 0, -1, 2)$ and $\{S(0, \varepsilon t, t, 0, \delta t, 0)\}_{t\in\mathcal H_1, \varepsilon, \delta=\pm1}$ up to isomorphism (see Proposition~\ref{pro:8.5}).  Also recall that all algebras in $(\mathcal{ECS}_{011})_3$ are classified into the family\\ $\{S(0, -t, t, t, -t, 0)\}_{ t\in\{-4\}\cup \mathcal H'_5}$ up to isomorphism (see Proposition~\ref{pro:9.1} (i)).
\vspace{2mm}

(A-1) We show that $S(0, 1, 1, 0, -1, 2)\ncong S(0, -t, t, t, -t, 0)$ for all $t\in K^*$.  In fact, let $t\in K^*$ and assume $S(0, 1, 1, 0, -1, 2)\cong S(0, -t, t, t, -t, 0)$.  Then we see from Lemma~\ref{lem:iso} that there are $x, y, z, w\in K$ with $\left|\begin{matrix}x&y\\z&w\end{matrix}\right|\ne0$ such that 
\begin{equation}\label{eq:1&3} % (27)
\left\{\begin{array}{@{\,}lll}
z=0, \\
x^2-ty^2+txy=w, \\
-tw^2=w, \\
txw=x, \\
-tyw+txw=y, \\
-txw=-x, \\
-tyw=-y+2w
\end{array} \right. 
\end{equation}
holds.  By the first equation of (\ref{eq:1&3}), we obtain $x\ne0$ and $w\ne0$.  Then by the third and 4th equations of (\ref{eq:1&3}), we obtain $w=-1/t$ and $w=1/t$, a contradiction.
\vspace{2mm}

(A-2) Let $t, t'\in K^*$ and $\varepsilon, \delta\in\{1, -1\}$.   Then we show that 
\[
S(0, \varepsilon t, t, 0, \delta t, 0)\cong S(0, -t', t', t', -t', 0)\Leftrightarrow \left\{\begin{array}{@{\,}lll}\varepsilon=-1, \\\delta=-1, \\t'\ne-4\end{array} \right. \, {\rm{and}}\,\, \frac{t}{t'(4+t')}\in(K^*)^2.
\]
In fact, suppose $S(0, \varepsilon t, t, 0, \delta t, 0)\cong S(0, -t', t', t', -t', 0)$.  By  Lemma~\ref{lem:iso}, there are $x, y, z, w\in K$ with $\left|\begin{matrix}x&y\\z&w\end{matrix}\right|\ne0$ such that 
\begin{equation}\label{eq:1&3another} % (28)
\left\{\begin{array}{@{\,}lll}
z=0, \cdots(28-1)\\
x^2-t'y^2+t'xy=w, \cdots(28-2)\\
-t'w=\varepsilon t, \cdots(28-3)\\
t'w=t, \cdots(28-4)\\
-t'yw+t'xw=ty, \cdots(28-5)\\
-t'w=\delta t, \cdots(28-6)\\
-t'yw=\delta ty. \cdots(28-7)
\end{array} \right. 
\end{equation}
By (28-3) and (28-4), we obtain $\varepsilon=-1$.  By (28-4) and (28-6), we obtain $\delta=-1$.  By (28-4) and (28-5), we obtain $x=2y$.  Then by (28-2) and (28-4), we obtain $t'(4+t')y^2=t$, hence $t'\ne-4$ and $\frac{t}{t'(4+t')}\in(K^*)^2$. 

Conversely suppose that $\left\{\begin{array}{@{\,}lll}\varepsilon=-1, \\\delta=-1, \\t'\ne-4\end{array} \right. $ and $\frac{t}{t'(4+t')}\in(K^*)^2$.  Then there is $\alpha\in K^*$ such that $\frac{t}{t'(4+t')}=\alpha^2$.  Define $x=2\alpha, y=\alpha, z=0$ and $w=t/t'$.
 Then $\left|\begin{matrix}x&y\\z&w\end{matrix}\right|=2\alpha t/t'\ne0$.  Moreover we can easily see that the scalars $x, y, z$ and $w$ satisfy (\ref{eq:1&3another}), and hence $S(0, \varepsilon t, t, 0, \delta t, 0)\cong S(0, -t', t', t', -t', 0)$ with a transformation matrix $\begin{pmatrix}2\alpha&\alpha\\0&t/t'\end{pmatrix}$.  
\vspace{2mm}

(A-3)  We show that $S(0, 4, -4, -4, 4, 0)\ncong S(0, \varepsilon t, t, 0, \delta t, 0)$ for all $t\in\mathcal H_1$ and  $\varepsilon, \delta=\pm1$.  To do this, suppose that  $S(0, 4, -4, -4, 4, 0)\cong S(0, \varepsilon t, t, 0, \delta t, 0)$ for some $t\in\mathcal H_1$ and $\varepsilon, \delta=\pm1$.  By (A-2)  (take $t'=-4$), we have $\varepsilon=\delta=-1$ and $-4\ne-4$, a contradiction.
\vspace{2mm}

Then we have the following.

\vspace{1mm}

\begin{pro}\label{pro:11.1} % Proposition 11.1
Suppose ${\rm{char}}\, K\ne2$. Then all algebras in $(\mathcal{ECS}_{011})_1\cup(\mathcal{ECS}_{011})_3$ are classified into 3 families $S(0, 1, 1, 0, -1, 2), \left\{S(0, \varepsilon t, t, 0, \delta t, 0)\right\}_{t\in\mathcal H_1, \varepsilon, \delta=\pm1}$ and $S(0, 4, -4, -4, 4, 0)$ up to isomorphism.
\end{pro}

\begin{proof}
By (A-1) and (A-3), we have only to show that $S(0, -t', t', t', -t', 0)$ is isomorphic to some algebra in $(\mathcal{ECS}_{011})_1$ whenever $t'\in\mathcal H'_5$.  To do this, let $t'\in\mathcal H'_5$ and put $t=4t'(4+t')$.  Then $t\in K^*$ and hence  we can take $t_1\in\mathcal H_1$ with $t_1\sim_1 t$ because $\mathcal H_1$ is a complete representative system of $K^*/\sim_1$.  Therefore we have
\[
\frac{t_1}{t'(4+t')}=\frac{t}{t'(4+t')}\cdot\frac{t_1}{t}=4\cdot\frac{t_1}{t}\in(K^*)^2\cdot(K^*)^2=(K^*)^2.
\] 
By (A-2), we see that $S(0, \varepsilon t_1, t_1, 0, \delta t_1, 0)\cong S(0, -t', t', t', -t', 0)$ for $\varepsilon=\delta=-1$.  Since $S(0, -t_1, t_1, 0, -t_1, 0)$ is an algebra in $(\mathcal{ECS}_{011})_1$, we have achieved what we aimed to show.
\end{proof}

(II) Relationship between $(\mathcal{ECS}_{011})_1$ and $(\mathcal{ECS}_{011})_4$
\vspace{2mm}

Recall that all algebras in $(\mathcal{ECS}_{011})_{1}$ are classified into 2 families $S(0, 1, 1, 0, -1, 2)$ and $\{S(0, \varepsilon t, t, 0, \delta t, 0)\}_{t\in\mathcal H_1, \varepsilon, \delta=\pm1}$ up to isomorphism.  Also recall that  all algebras in $(\mathcal{ECS}_{011})_4$ are classified into the family $\{S(0, (1+t)^2/4, (t^2-1)/4, 1, (1-t^2)/4, t)\}_{t\in K^* \backslash \{\pm1\}}$ up to isomorphism (see Proposition~\ref{pro:10.1}). 
\vspace{2mm}

(B-1) We show that 
\[
S(0, 1, 1, 0, -1, 2)\ncong S(0, (1+t')^2/4, (t'^2-1)/4, 1, (1-t'^2)/4, t')
\]
for all $t'\in K^* \backslash \{\pm1\}$.  In fact, otherwise, we see from Lemma~\ref{lem:iso} that there are $x, y, z, w\in K$ with $\left|\begin{matrix}x&y\\z&w\end{matrix}\right|\ne0$ such that  
\begin{equation}\label{eq:1&4} % (29)
\left\{\begin{array}{@{\,}lll}
z=0, \\
x^2+(1+t')^2y^2/4+(1+t')xy=w, \\
(1+t')^2w/4=1, \\
(t'^2-1)w/4=1, \\
(1+t')^2yw/4+xw=y, \\
(1-t'^2)w/4=-1, \\
y+t'xw=-y+2w
\end{array} \right. 
\end{equation}
holds for some $t'\in K^* \backslash \{\pm1\}$.  By the third and 4th equations of (\ref{eq:1&4}), we obtain $(1+t')^2=t'^2-1$, hence $t'=-1$, so $0=-1$ by the 6th equation of (\ref{eq:1&4}), a contradiction.
\vspace{2mm}

(B-2) We show that
\[
S(0, \varepsilon t, t, 0, \delta t, 0)\ncong S(0, (1+t')^2/4, (t'^2-1)/4, 1, (1-t'^2)/4, t')
\]
for all $t\in\mathcal H_1, \varepsilon, \delta=\pm1$ and $t'\in K^* \backslash \{\pm1\}$.  In fact, otherwise, we see from Lemma~\ref{lem:iso}  that there are $x, y, z, w\in K$ with $\left|\begin{matrix}x&y\\z&w\end{matrix}\right|\ne0$ such that  
\begin{equation}\label{eq:1&4another} % (30)
\left\{\begin{array}{@{\,}lll}
z=0, \\
x^2+q'y^2+(1+t')xy=w, \\
q'w=\varepsilon t, \\
a'w=t, \\
q'yw+xw=ty, \\
c'w=\delta t, \\
q'yw+t'xw=\delta ty
\end{array} \right. 
\end{equation}
holds for some $\varepsilon, \delta=\pm1$ and $t, t'\in K^*$ with $t'\ne\pm1$, where $q'=(1+t')^2/4, a'=(t'^2-1)/4$ and $c'=(1-t'^2)/4$.  By the third and 4th equations of (\ref{eq:1&4another}), we obtain $\varepsilon q'=a'$, that is, $\varepsilon(1+t')^2=t'^2-1$, hence $\varepsilon(1+t')=t'-1$, a contradiction because $t'\ne0$ and $\varepsilon =\pm1$.

By (B-1) and (B-2), we have the following.

\begin{pro}\label{pro:11.2} % Proposition 11.2
Suppose ${\rm{char}}\, K\ne2$.  Then all algebras in $(\mathcal{ECS}_{011})_1\cup(\mathcal{ECS}_{011})_4$ are classified into 3 families $S(0, 1, 1, 0, -1, 2), \{S(0, \varepsilon t, t, 0, \delta t, 0)\}_{t\in\mathcal H_1, \varepsilon, \delta=\pm1}$ and $\{S(0, (1+t)^2/4, (t^2-1)/4, 1, (1-t^2)/4, t)\}_{t\in K^* \backslash \{\pm1\}}$ up to isomorphism.
\end{pro}
\vspace{2mm}

(III) Relationship between $(\mathcal{ECS}_{011})_3$ and $(\mathcal{ECS}_{011})_4$
\vspace{2mm}

Recall that all algebras in $(\mathcal{ECS}_{011})_3$ are classified into the family\\ $\{S(0, -t, t, t, -t, 0)\}_{ t\in\{-4\}\cup \mathcal H'_5}$ up to isomorphism.   Also recall that  all algebras in $(\mathcal{ECS}_{011})_4$ are classified into the family $\{S(0, (1+t)^2/4, (t^2-1)/4, 1, (1-t^2)/4, t)\}_{t\in K^* \backslash \{\pm1\}}$ up to isomorphism. 
 We show that 
\[
S(0, -t, t, t, -t, 0)\ncong S(0, (1+t')^2/4, (t'^2-1)/4, 1, (1-t'^2)/4, t')
\]
for all $t\in\{-4\}\cup \mathcal H'_5$ and $t'\in K^* \backslash \{\pm1\}$.  In fact, otherwise, we see from Lemma~\ref{lem:iso} that there are $x, y, z, w\in K$ with $\left|\begin{matrix}x&y\\z&w\end{matrix}\right|\ne0$ such that 
\begin{equation}\label{eq:3&4} % (31)
\left\{\begin{array}{@{\,}lll}
z=0, \\
x^2+(1+t')^2y^2/4+(1+t')xy=w, \\
(1+t')^2w/4=-t, \\
(t'^2-1)w/4=t, \\
(1+t')^2yw/4+xw=ty+tw, \\
(1+t')^2yw/4+t'xw=-ty
\end{array} \right. 
\end{equation}
holds for some $t\in\{-4\}\cup \mathcal H'_5$ and $t'\in K^* \backslash \{\pm1\}$.  By the third and 4th equations of (\ref{eq:3&4}), we obtain $(1+t')^2=1-t'^2$, hence $t'=-1$, a contradiction.  Therefore we have the following.
 
\begin{pro}\label{pro:11.3} % Proposition 11.3
Suppose ${\rm{char}}\, K\ne2$.  Then all algebras in $(\mathcal{ECS}_{011})_3\, \cup(\mathcal{ECS}_{011})_4$ are classified into 2 families $\{S(0, -t, t, t, -t, 0)\}_{ t\in\{-4\}\cup \mathcal H'_5}$ and\\ $\{S(0, (1+t)^2/4, (t^2-1)/4, 1, (1-t^2)/4, t)\}_{t\in K^* \backslash \{\pm1\}}$ up to isomorphism.
\end{pro}
\vspace{1mm}

By leveraging Propositions 11.1, 11.2 and 11.3, we can easily prove
Theorem~\ref{thm:char-not=2}.  Indeed, put
\[
\mathcal A=\{S(0, 1, 1, 0, -1, 2), S(0, 4, -4, -4, 4, 0\},  \mathcal B=\{S(0, \varepsilon t, t, 0, \delta t, 0) : t \in\mathcal H_1, \varepsilon, \delta=\pm1\}
\]
and
\[
 \mathcal C=\{S(0, (1+t)^2/4, (t^2-1)/4, 1, (1-t^2)/4, t) : t\in K^* \backslash \{\pm1\}\}.
\]
Then $\mathcal A\cup\mathcal B\cup\mathcal C\subseteq\mathcal{ECS}_{011}$ and any two algebras in $\mathcal A\cup\mathcal B\cup\mathcal C$ are not isomorphic to each other from Propositions 11.1, 11.2 and 11.3.  Furthermore, any algebra in $\mathcal{ECS}_{011}$ is isomorphic to some algebra in $\mathcal A\cup\mathcal B\cup\mathcal C$ from the same propositions.

\section{Proof of Theorem~\ref{thm:char=2}}\label{sec:proof2} % section 12

In this section, we assume that ${\rm{char}}\, K=2$.  
\vspace{2mm}

(I) Relationship between $(\mathcal{ECS}_{011})_1$ and $(\mathcal{ECS}_{011})_3$
\vspace{2mm}

Recall that all algebras in $(\mathcal{ECS}_{011})_{1}$ are classified into 2 families $\{S(0, t, t, 0, t, 1)\}_{t\in\mathcal H_2}$ and $\{S(0, t, t, 0, t, 0)\}_{t\in\mathcal H_3}$ up to isomorphism (see Proposition~\ref{pro:8.6}).  Also recall that all algebras in $(\mathcal{ECS}_{011})_3$ are classified into the family $\{S(0, t, t, t, t, 0)\}_{ t\in{\mathcal H_4}}$ up to isomorphism
(see Proposition~\ref{pro:9.1} (ii)).

We assert that $S(0, t, t, 0, t, \delta)\ncong S(0, t', t', t', t', 0)$ for all $t, t'\in K^*$ and $\delta\in\{0, 1\}$.  In fact, let $t, t'\in K^*$ and $\delta\in\{0, 1\}$.  We assume $S(0, t, t, 0, t, \delta)\cong S(0, t', t', t', t', 0)$ and derive a contradiction as follows.  By Lemma~\ref{lem:iso}, there are $x, y, z, w\in K$ with $\left|\begin{matrix}x&y\\z&w\end{matrix}\right|\ne0$ such that  
\begin{equation}\label{eq:1&3-char2} % (32) 
\left\{\begin{array}{@{\,}lll}
z=0, \\
x^2+t'y^2+t'xy=w, \\
t'w^2=tw, \\
t'xw=tx, \\
t'yw+t'xw=ty, \\
t'xw=tx, \\
t'yw=ty+\delta w
\end{array} \right. 
\end{equation}
holds.  
 By the first equation of (\ref{eq:1&3-char2}), we obtain $x\ne0$ and $w\ne0$.  Then by the third equation of (\ref{eq:1&3-char2}), we obtain $t'w=t$, hence $ty+tx=ty$ by the 5th equation of (\ref{eq:1&3-char2}).  Then we have $tx=0$, a contradiction.  Therefore, we have the following.
\vspace{2mm}

\begin{pro}\label{pro:12.1} % Proposition 12.1
Suppose ${\rm{char}}\, K=2$.  Then all algebras in $(\mathcal{ECS}_{011})_1\cup(\mathcal{ECS}_{011})_3$ are classified into 3 families $\{S(0, t, t, 0, t, 1)\}_{t\in\mathcal H_2}, \{S(0, t, t, 0, t, 0)\}_{t\in\mathcal H_3}$ and \\$\{S(0, t, t, t, t, 0)\}_{ t\in{\mathcal H_4}}$ up to isomorphism.
\end{pro}
\vspace{2mm}

(II) Relationship between $(\mathcal{ECS}_{011})_1$ and $(\mathcal{ECS}_{011})_4$
\vspace{2mm}

Recall that all algebras in $(\mathcal{ECS}_{011})_{1}$ are classified into 2 families $\{S(0, t, t, 0, t, 1)\}_{t\in\mathcal H_2}$ and $\{S(0, t, t, 0, t, 0)\}_{t\in\mathcal H_3}$ up to isomorphism.  Also recall that all algebras in $(\mathcal{ECS}_{011})_4$ are classified into the family $\left\{S\left(0, \frac{t^2}{1+t^2}, \frac{t}{1+t^2}, 1, \frac{t}{1+t^2}, 1\right)\right\}_{t\in K^*\backslash \{1\}}$ up to isomorphism (see Proposition~\ref{pro:10.2}).

We assert that $S(0, t, t, 0, t, \delta)\ncong S\left(0, \frac{t'^2}{1+t'^2}, \frac{t'}{1+t'^2}, 1, \frac{t'}{1+t'^2}, 1\right)$ holds for all $t\in K^*, \delta\in\{0, 1\}$ and $t'\in K^*\backslash \{1\}$.  In fact, otherwise, we see from Lemma~\ref{lem:iso} that there are $x, y, z, w\in K$ with $\left|\begin{matrix}x&y\\z&w\end{matrix}\right|\ne0$ such that  
\begin{equation}\label{eq:1&4-char2} % (33)
\left\{\begin{array}{@{\,}lll}
z=0, \\
x^2+\frac{t'^2}{1+t'^2}y^2=w, \\
\frac{t'^2}{1+t'^2}w=t, \\
\frac{t'}{1+t'^2}w=t, \\
\frac{t'^2}{1+t'^2}yw+xw=ty, \\
\frac{t'^2}{1+t'^2}yw+xw=ty+\delta w
\end{array} \right. 
\end{equation}
holds for some $t\in K^*, \delta\in\{0, 1\}$ and $t'\in K^*\backslash \{1\}$.   By the third and 4th equations of (\ref{eq:1&4-char2}), we obtain  $t'^2=t'$, that is, $t'=1$, a contradiction.  

By our assertion, we have the following.

\begin{pro}\label{pro:12.2} % Proposition 12.2
Suppose ${\rm{char}}\, K=2$.  Then all algebras in $(\mathcal{ECS}_{011})_1\cup(\mathcal{ECS}_{011})_4$ are classified into 3 families $\{S(0, t, t, 0, t, 1)\}_{t\in\mathcal H_2}, \{S(0, t, t, 0, t, 0)\}_{t\in\mathcal H_3}$ and\\
 $\left\{S\left(0, \frac{t^2}{1+t^2}, \frac{t}{1+t^2}, 1, \frac{t}{1+t^2}, 1\right)\right\}_{t\in K^*\backslash \{1\}}$ up to isomorphism.
\end{pro}
\vspace{2mm}

(III) Relationship between $(\mathcal{ECS}_{011})_3$ and $(\mathcal{ECS}_{011})_4$
\vspace{2mm}

Recall that all algebras in $(\mathcal{ECS}_{011})_3$ are classified into the family $\{S(0, t, t, t, t, 0)\}_{ t\in {\mathcal H_4}}$ up to isomorphism.   Also recall that all algebras in $(\mathcal{ECS}_{011})_4$ are classified into the family $\left\{S\left(0, \frac{t^2}{1+t^2}, \frac{t}{1+t^2}, 1, \frac{t}{1+t^2}, 1\right)\right\}_{t\in K^*\backslash \{1\}}$ up to isomorphism.  We show that $S(0, t, t, t, t, 0)\ncong S\left(0, \frac{t'^2}{1+t'^2}, \frac{t'}{1+t'^2}, 1, \frac{t'}{1+t'^2}, 1\right)$ holds for all $t\in {\mathcal H_4}$ and $t'\in K^*\backslash \{1\}$.  In fact, otherwise, we see from Lemma~\ref{lem:iso}  that there are $x, y, z, w\in K$ with $\left|\begin{matrix}x&y\\z&w\end{matrix}\right|\ne0$ such that  
\begin{equation}\label{eq:3&4-char2} % (34)
\left\{\begin{array}{@{\,}lll}
z=0, \\
x^2+\frac{t'^2}{1+t'^2}y^2=w, \\
\frac{t'^2}{1+t'^2}w=t, \\
\frac{t'}{1+t'^2}w=t, \\
\frac{t'^2}{1+t'^2}yw+xw=ty+tw, \\
\frac{t'^2}{1+t'^2}yw+xw=ty
\end{array} \right. 
\end{equation}
holds for some $t\in{\mathcal H_4}$ and $t'\in K^*\backslash \{1\}$.   By the third and 4th equations of (\ref{eq:3&4-char2}), we obtain $t'^2=t'$, hence $t'=1$, a contradiction.  Therefore, we have the following.

\begin{pro}\label{pro:12.3} % Proposition 12.3
Suppose ${\rm{char}}\, K=2$.  Then all algebras in $(\mathcal{ECS}_{011})_3\cup(\mathcal{ECS}_{011})_4$ are classified into 2 families $\{S(0, t, t, t, t, 0)\}_{ t\in {\mathcal H_4}}$ and\\ $\left\{S\left(0, \frac{t^2}{1+t^2}, \frac{t}{1+t^2}, 1, \frac{t}{1+t^2}, 1\right)\right\}_{t\in K^*\backslash \{1\}}$ up to isomorphism.
\end{pro}

\vspace{1mm}

By leveraging Propositions~\ref{pro:12.1}, \ref{pro:12.2} and \ref{pro:12.3}, we can easily prove Theorem~\ref{thm:char=2}.  Indeed, put $\mathcal A=\{S(0, t, t, 0, t, 1) : t\in\mathcal H_2\}$, $\mathcal B=\{S(0, t, t, 0, t, 0) : t\in\mathcal H_3\}$,  $\mathcal C=\{S(0, t, t, t, t, 0) : t\in {\mathcal H_4}\}$ and $ \mathcal D=\left\{S\left(0, \frac{t^2}{1+t^2}, \frac{t}{1+t^2}, 1, \frac{t}{1+t^2}, 1\right): t\in K^*\backslash \{1\}\right\}$.
Then $\mathcal A\cup\mathcal B\cup\mathcal C\cup\mathcal D\subseteq\mathcal{ECS}_{011}$ and any two algebras in $\mathcal A\cup\mathcal B\cup\mathcal C\cup\mathcal D$ are not isomorphic to each other from Propositions~\ref{pro:12.1}, \ref{pro:12.2} and \ref{pro:12.3}.  Furthermore, any algebra in $\mathcal{ECS}_{011}$ is isomorphic to some algebra in $\mathcal A\cup\mathcal B\cup\mathcal C\cup\mathcal D$ from the same propositions.

\end{document}